\pdfoutput=1
\documentclass[11pt,a4paper]{article}
\usepackage[margin=1in]{geometry}
\usepackage[utf8]{inputenc}
\usepackage[T1]{fontenc}
\usepackage{lmodern}
\usepackage[english]{babel}
\usepackage{amsmath,amssymb,amsthm,mathtools}
\usepackage{bm}
\usepackage{enumitem}
\usepackage{microtype}
\allowdisplaybreaks[1]
\usepackage{booktabs}
\usepackage{multirow}
\usepackage{makecell}
\usepackage{mathrsfs}
\usepackage{graphicx}
\usepackage{caption}
\usepackage{subcaption}
\usepackage{placeins}
\usepackage{latexsym,amsfonts,amsbsy}
\usepackage{algorithmic}
\usepackage{float}
\usepackage[colorlinks=true,linkcolor=blue,citecolor=blue,urlcolor=blue]{hyperref}

\newtheorem{definition}{Definition}[section]
\newtheorem{assumption}{Assumption}
\newtheorem{lemma}{Lemma}
\newtheorem{proposition}{Proposition}
\newtheorem{theorem}{Theorem}

\theoremstyle{remark}
\newtheorem{remark}{Remark}
\newtheorem{example}{\textbf{\textsl{Example}}}
\newcounter{algorithm}
\newcommand{\R}{\mathbb{R}}

\newcommand{\Nzero}{\mathbb{N}_0}
\newcommand{\dq}{\,\mathrm{d}_q}
\newcommand{\CD}[1]{{}^{\mathrm C}D_q^{#1}}
\newcommand{\Jq}[1]{J_q^{#1}}
\newcommand{\norm}[1]{\left\lVert #1\right\rVert}
\newcommand{\abs}[1]{\left\lvert #1\right\rvert}

\numberwithin{equation}{section}

\title{A singularity detecting uniformly convergent quantum wavelet framework for weakly singular nonlocal Emden--Fowler systems with quantum Mittag--Leffler sources}
\author{%
Arushi\textsuperscript{1}, Pranay Goswami\textsuperscript{2}, Pratibhamoy Das\textsuperscript{3,*}, Sudarshan Santra\textsuperscript{4}\\[8pt]
\textsuperscript{1}\textit{Sri Venkateswara College, University of Delhi, Delhi, India}\\
\texttt{arushi.garvita@gmail.com}, \texttt{arushi@svc.ac.in}\\[3pt]
\textsuperscript{2}\textit{Dr.\ B.\ R.\ Ambedkar University Delhi, Delhi, India}\\
\texttt{pranay@aud.ac.in}, \texttt{pranaygoswami83@gmail.com}\\[3pt]
\textsuperscript{3}\textit{Department of Mathematics, Indian Institute of Technology Patna, Bihar, India}\\
\texttt{pratibhamoy@gmail.com}\\[3pt]
\textsuperscript{4}\textit{Department of Mathematics, Birla Institute of Technology and Science Pilani, Rajasthan, India}\\
\texttt{sudarshan.santra@pilani.bits-pilani.ac.in}\\[6pt]
\textsuperscript{*}Corresponding author%
}
\date{}
\begin{document}
\maketitle
\begin{abstract}
\noindent We develop a quantum wavelet-based computational framework for coupled systems of Caputo-type quantum ($q$-) fractional Emden--Fowler problems with $q$-Mittag--Leffler source terms, where the  solutions exhibit singular behavior in time. Here $q\in(0,1)$ denotes the parameter of $q$-calculus. The highest-order Caputo $q$-fractional derivatives are expanded through the Haar wavelet basis, and the lower-order fractional derivatives are approximated through Jackson $q$-fractional integrals of Haar wavelets. The continuous problem reduces to a nonlinear algebraic system on the midpoint collocation points. Under suitable structural assumptions, the Banach fixed-point theorem yields the existence and uniqueness of the continuous solutions and wellposedness of the discrete solutions on collocation points. The nonlinearity of the given problem is dealt with by the damped Newton's iteration. In addition, We prove stability of the scheme and derive rigorous convergence estimates in a weighted norm that remains uniformly controlled despite the appearance of singular weight at $t=0$ in the continuous model. When the highest-order fractional derivatives are H\"older continuous with exponent $\sigma$, the mathematical error of weakly singular solution and its lower-order fractional derivatives become $\mathcal{O}(\tau_L^{\sigma})$, where $\tau_L=2^{-(L+1)}$ denotes the width of the finest subintervals for maximum level of resolution $L$. First-order convergence is guaranteed under stronger regularity with $\sigma=1$. Several test cases are considered to confirm the accuracy and robustness of the present approach, which strongly agrees with the theoretical estimates. In particular, first-order convergence is observed even where the stronger regularity assumption is not satisfied.
\end{abstract}

\section{Introduction}
Quantum calculus, also known as $q$-calculus, generalizes classical calculus through a framework that does not explicitly rely on the conventional notion of limits. In recent years, it has been proven to be an active area of mathematical research, encompassing a broad range of theoretical foundations and applications \cite{KacCheung2002,Ernst2012,Ahmad2016}. An important extension of this framework is $q$-fractional calculus, which introduces $q$-analogues of classical fractional integrals and derivatives \cite{annaby2012fractional}, drawing on the rich theory of $q$-special functions and basic hypergeometric series \cite{GasperRahman2004,Ismail2005}. Similar to classical fractional calculus \cite{Kilbas2006,Diethelm2010}, $q$-fractional calculus incorporates nonlocal operators and associated with weakly singular kernels. In the classical framework, kernels of this type commonly display a power-law behaviour characterized by the term  $(t-s)^{\alpha-1}$, where the parameter $\alpha$ lies in $(0,1)$. The inherent memory effects encoded in fractional operators, combined with the weak regularity of solutions near $t=0$, give rise to considerable difficulties from both analytical and computational perspectives. Moreover, having singularities in the solution's profiles of fractional differential and integro-differential equations \cite{SantraAnalysis2025}, leads to a low order of accuracy of conventional high-order numerical schemes that rely on high regularity assumptions. This provides a strong motivation for exploring $q$-fractional extensions of singular differential equations that emerge in mathematical physics, with particular attention to generalizations of the Lane–Emden and Emden–Fowler equations. Their classical forms, which involve singular behavior at the origin, have been extensively investigated in the context of polytropic stellar models \cite{KumarGoswami2023,Shahni2021,Singh2021}.

The numerical solution of nonlocal or integro-differential equations has traditionally relied on finite difference and finite element methods. Finite difference schemes, such as the L1 \cite{SANTRA202387} and L2-1$\sigma$ \cite{SantraEnhancing2025} approximations of the Caputo derivative, are simple to implement and achieve high accuracy when the solution is sufficiently smooth \cite{Alikhanov2021}. Finite element methods offer greater geometric flexibility and a well-developed error analysis in variational settings \cite{SANTRA2025}. Both families, however, face serious limitations in the fractional setting. The nonlocal, weakly singular memory kernel forces each time step to depend on the entire solution history, which makes the computational cost and storage grow rapidly with the number of time levels. Moreover, typically the solution derivatives become unbounded near the initial time, so the optimal convergence rates are lost on uniform meshes unless graded meshes or other corrective strategies are used \cite{SantraEnhancing2025}. Adaptive mesh strategy requires strong stability analysis and can deteriorate the condition number of the resulting algebraic system.

Local irregularities are often captured by orthogonal wavelet functions having compact support. In addition the multiresolution structure of the wavelet approximations often reduce the computational costs of fractional derivatives as well as nonlocal operators like integrals \cite{Ray2012}. Among them, the Haar wavelet collocation method (HWCM) is particularly attractive because of its simplicity, its fast computation, and the availability of explicit integration matrices. It has been applied successfully to weakly singular initial-boundary value problems \cite{Dassantra25} and single Emden--Fowler equations \cite{KumarGoswami2023}. Its convergence theory also requires less regularity of the solution than classical schemes do \cite{Majak2015}, which is especially useful for problems with weak singularities. Recently, \cite{kumar2024numerical} extended HWCM to coupled Emden--Fowler systems and compared the resulting approximations with artificial neural-network approaches; related neural-network investigations may also be found in \cite{PATAWARI2027109554,SantrafI-PIDNNs2027}. However, in this work, we focus on a $q$-analogue of the Haar wavelet collocation method that is designed to handle the weak singularities arising in $q$-fractional problems. By exploiting the low regularity requirements of Haar wavelets, the proposed approach aims to retain accuracy for solutions that are non-smooth near the initial time, without relying on graded meshes.    
	
The works on $q$-fractional operators and $q$-Mittag-Leffler functions \cite{Stankovic2009,rajkovic2007qanalogues,AbdeljawadBaleanu2011} provide the theoretical base for the source terms employed in this paper, confirming that $q$-special functions remain an important area of foundational development. The $q$-fractional operators and $q$-Mittag--Leffler functions require quantum calculus which is an advancement of usual fractional operators, like Caputo fractional derivative. Our research is motivated on the use of simple Haar wavelet collocation for the coupled $q$-fractional Emden--Fowler-type nonlinear systems where the time singularity appears and to provide a uniformly convergent method for time singular problems through quantum calculus. One can also look into the recent works on $q$-calculus \cite{IntegralRepresentations2026, HALBURD2026130336, An25} for different classes of fractional differential equations.

Despite these developments, comparatively few studies are available on $q$-fractional differential equations. For instance, Zhang and Guo \cite{zhang2020solution} established the existence and stability results for $q$-fractional differential equations, providing the theoretical foundation for numerical methods in this domain. However, the present work deals with the quantum wavelet framework for weakly singular nonlocal Emden--Fowler systems with quantum Mittag-Leffler sources for which the solutions have unbounded time derivatives at $t=0$. To the best of our knowledge, this is the first attempt to treat such singular $q$-fractional systems within a wavelet setting.
	
In this paper, we develop a Haar wavelet collocation method for a coupled system of Caputo-type $q$-fractional Emden--Fowler equations with $q$-Mittag--Leffler source terms. The highest-order $q$-fractional derivatives are approximated using Haar wavelets, while the state variables and lower-order derivatives are reconstructed through Jackson $q$-fractional integration. The resulting nonlinear algebraic system is solved numerically, and existence, uniqueness, stability, and convergence of the method are established. Numerical examples are presented to illustrate the accuracy and first-order convergence behavior of the proposed approach. The main contribution of the present work is summarized as follows:
    \begin{itemize}
        \item Develop a novel quantum wavelet-based collocation framework for coupled Caputo-type $q$-fractional Emden--Fowler systems with $q$-Mittag--Leffler sources having solutions with weak temporal singularities.
        
        \item Using Jackson $q$-fractional integration of Haar basis functions and midpoint collocation, we reduce the coupled system to a finite-dimensional nonlinear algebraic problem and solve it via a backtracking damped Newton method.

        \item Establish the existence, uniqueness, and stability of the continuous solution and the discrete solution on collocation points using the Banach fixed-point theorem.
        
        \item Derive rigorous convergence estimates of order $\mathcal{O}(\tau_L^\sigma)$ for $0<\sigma\leq1$, yielding first-order convergence when $\sigma=1$, with uniform control of the weighted error near $t=0$.
        
        \item Validate the proposed method with several test cases, demonstrating its accuracy and robustness, with first-order convergence which are observed even for cases where Lipschitz regularity is not satisfied.
    \end{itemize}
The remainder of the paper is organized as follows. Section~\ref{Section_Preleminaries} reviews the essential concepts of $q$-calculus, $q$-fractional calculus, $q$-Mittag--Leffler functions, and Haar wavelets. Section~\ref{Section_ProblemFormulation} formulates the coupled $q$-fractional Emden--Fowler system and derives its equivalent Volterra-type $q$-integral representation. Section~\ref{Section_Methodology} presents the proposed quantum calculus-adapted wavelet discretization, leading to the midpoint collocation scheme and the resulting nonlinear algebraic system. Section~\ref{Section_Qualitative} establishes the existence and uniqueness of the continuous solution and the discrete collocation solution. We develop the stability and convergence analysis in Section~\ref{Section_Convergence}, accounting for the weak temporal singularity at $t=0$, and establishes an error estimate of order $\mathcal{O}(\tau_L^\sigma)$ under H\"{o}lder regularity of exponent $\sigma\in(0,1]$, with first-order convergence guaranteed when $\sigma=1$. The damped Newton's iteration scheme to deal with the nonlinearity of the system has been discussed in Section \ref{Sec_Newton}. Finally, Section~\ref{Section_Numerical_Examples} presents a couple of examples to assess the accuracy and robustness of the proposed method and examine its observed convergence behavior, while Section~\ref{Section_Conclusion} summarizes the main findings and outlines possible directions for future research.

\section{Preliminaries}\label{Section_Preleminaries}
Throughout the paper, $0<q<1,~ I=[0,1],$
and $\mathbb{R}$ denotes the set of real numbers. The independent variable is denoted by $t\in I$. We first recall some fundamental concepts and definitions from $q$-calculus that will be used in the subsequent analysis. 
\begin{definition}[q-derivative \cite{KacCheung2002}]
		For $f:[0,1]\to\R$ and $t\in(0,1]$, define
		$
			D_qf(t)=\frac{f(t)-f(qt)}{(1-q)t}.
		$
		Set $D_qf(0)=\lim_{t\to0^+}D_qf(t)$. Higher-order $q$-derivatives are defined recursively by $D_q^0f=f$ and $D_q^nf=D_q(D_q^{n-1}f)$.
\end{definition}
	
\begin{definition}[Jackson $q$-integral \cite{KacCheung2002,Ernst2012}]\label{Def_q_Integral}
		For $b\in[0,1]$, provided that the series converges absolutely,
		$
			\int_0^b f(s)\dq s=(1-q)b\sum_{n=0}^{\infty}q^n f(bq^n).
		$
		For $0\le a<b\le1$, define
		$
			\int_a^b f(s)\dq s=\int_0^b f(s)\dq s-\int_0^a f(s)\dq s.
		$
\end{definition}
	
\begin{definition}[\cite{Stankovic2009,GasperRahman2004}]
		For $n\in\Nzero$ and $a\in\mathbb C$, let $(a;q)_0=1, (a;q)_n=\prod_{r=0}^{n-1}(1-aq^r)$ and $
			(a;q)_\infty=\prod_{r=0}^{\infty}(1-aq^r),$ whenever the infinite product converges. For $x>0$, $0\le y<x$, and $\nu\in\R$, define
		$(x-y)_q^\nu=x^\nu\frac{(y/x;q)_\infty}{(q^\nu y/x;q)_\infty},
			\label{eq:qpower}$
		provided that the denominator is nonzero.
\end{definition}
	
\begin{definition}[$q$-Gamma function \cite{KacCheung2002,GasperRahman2004,Ismail2005}]
		For $x>0$,
		$
			\Gamma_q(x)=(1-q)^{1-x}\frac{(q;q)_\infty}{(q^x;q)_\infty}.
		$
		It satisfies
		$\Gamma_q(x+1)=[x]_q\Gamma_q(x),[x]_q=\frac{1-q^x}{1-q}$ and $ \Gamma_q(1)=1.$
\end{definition}

\subsection{q-fractional Operators}
We now introduce the $q$-fractional operators required for the formulation and analysis of the proposed method.
\begin{definition}[Left-sided $q$-fractional integral \cite{Stankovic2009,annaby2012fractional}]
		Let $\alpha>0$. Whenever the Jackson integral exists, define
		$
			\Jq{\alpha}f(t)=\frac{1}{\Gamma_q(\alpha)}\int_0^t (t-qs)_q^{\alpha-1}f(s)\dq s,  t\in(0,1].
			\label{eq:qfracint}
		$
		Set $\Jq{\alpha}f(0)=0$ for $\alpha>0$, and $\Jq{0}f=f$.
\end{definition}
	
\begin{definition}[Caputo-type $q$-fractional derivative \cite{Stankovic2009,rajkovic2007qanalogues}]
		Let $\alpha>0$ and $n=\lceil\alpha\rceil$. The Caputo-type $q$-fractional derivative is
		$
			\CD{\alpha}f(t)=\Jq{n-\alpha}D_q^nf(t).
		$
		In particular, for $0<\alpha<1$,
		$
			\CD{\alpha}f(t)=\Jq{1-\alpha}D_qf(t)
			=\frac{1}{\Gamma_q(1-\alpha)}\int_0^t(t-qs)_q^{-\alpha}D_qf(s)\dq s.
		$
\end{definition}
	
\begin{lemma}[{\bf Composition identities} \cite{Stankovic2009,annaby2012fractional}]
		\label{lem:composition}
		Let the required $q$-derivatives and Jackson integrals exist, and assume absolute convergence of the series occurring in the compositions. For $\mu,\nu>0$,
		\begin{equation*}
			\Jq{\mu}\Jq{\nu}f=\Jq{\mu+\nu}f.
		\end{equation*}
		For $0<\alpha<1$,
		\begin{equation}
			\Jq{\alpha}\CD{\alpha}f=f-f(0),
			\qquad
			\CD{\alpha}\Jq{\alpha}f=f.
			\label{eq:composition}
		\end{equation}
		Consequently, if $0<\beta<\alpha<1$, then
		\begin{equation}
			\Jq{\alpha-\beta}\CD{\alpha}f=\CD{\beta}f.
			\label{eq:lowercomposition}
		\end{equation}
\end{lemma}
	
\subsection[\texorpdfstring{$q$-Mittag--Leffler Functions}{q-Mittag--Leffler Functions}]{$q$-Mittag--Leffler Functions}
The \(q\)-Mittag–Leffler functions provide a natural \(q\)-analogue of the classical Mittag–Leffler functions and play an important role in representing solutions of \(q\)-fractional differential equations. Moreover, as \(q\to1^{-}\), the \(q\)-Mittag–Leffler functions recover their classical counterparts, establishing a direct connection between quantum calculus and the standard fractional-calculus framework.	
\begin{definition}[$q$-Mittag--Leffler function \cite{AbdeljawadBaleanu2011,rajkovic2007qanalogues}]
		For $\alpha,\beta>0$, define
		\begin{equation}
			E_{q,\alpha,\beta}(z)=\sum_{k=0}^{\infty}\frac{z^k}{\Gamma_q(\alpha k+\beta)},
			\qquad \abs{z}<(1-q)^{-\alpha}.
			\label{eq:qML}
		\end{equation}
		When $\beta=1$, write $E_{q,\alpha}=E_{q,\alpha,1}$.
\end{definition}
	
\begin{proposition}
		\label{prop:qMLdisk}
		The series in \eqref{eq:qML} has radius of convergence $(1-q)^{-\alpha}$. On every compact subset of this disk, the series and each termwise derivative converge uniformly. Consequently, $E_{q,\alpha,\beta}$ is analytic in the disk, and $E_{q,\alpha}$ is Lipschitz continuous on every compact real interval contained in it.
\end{proposition}
	
\begin{proof}
		From the product representation of $\Gamma_q$,
		$
			\frac{\Gamma_q(\alpha k+\beta)}{\Gamma_q(\alpha(k+1)+\beta)}
			=(1-q)^{\alpha}
			\frac{(q^{\alpha(k+1)+\beta};q)_\infty}
			{(q^{\alpha k+\beta};q)_\infty}
			\longrightarrow (1-q)^{\alpha}.
	$
		The ratio test therefore gives the radius $(1-q)^{-\alpha}$. Uniform convergence of the power series and its derivatives on compact subsets is standard, and the Lipschitz assertion follows from boundedness of the derivative on a compact interval.
\end{proof}

\subsection{The Haar Family}
We next recall the Haar wavelet family and its basic properties, which form the foundation of the proposed numerical framework.
\begin{definition}[Haar functions]
		Let $h_1(t)=1$ for $0\le t<1$. For each level $\ell\in\Nzero$ and translation $k=1,\dots,2^\ell$, set
		$j=2^\ell+k,\xi_1=\frac{k-1}{2^\ell}, \xi_2=\frac{k-\tfrac12}{2^\ell}$ and $\xi_3=\frac{k}{2^\ell}$. Define
		\begin{equation*}
			h_j(t)=
			\begin{cases}
				1, & \xi_1\le t<\xi_2,\\
				-1, & \xi_2\le t<\xi_3,\\
				0, & \text{otherwise}.
			\end{cases}
		\end{equation*}
		The functions are extended to $t=1$ by their left-hand limits.
\end{definition}
\begin{remark}
The above family is orthogonal, but not normalized in $L^2(0,1)$ as $\norm{h_1}_{L^2(0,1)}^2=1$ and $\norm{h_j}_{L^2(0,1)}^2=2^{-\ell}$ for a wavelet with level of resolution $\ell$. For a maximum resolution level $L\ge0$, let $N_L=2^{L+1}$ and $\tau_L=\frac{1}{N_L}=2^{-(L+1)}$. Then $\{h_1,\dots,h_{N_L}\}$ span the space $V_L$ of functions that are constant on each dyadic cell $I_r=[(r-1)\tau_L,r\tau_L), \,r=1,\dots,N_L-1,$
	with the last cell $I_r=[(N_L-1)\tau_L,1]$.
\end{remark}

\section{Mathematical Model and \texorpdfstring{$q$-Integral}{q-Integral} Representation}\label{Section_ProblemFormulation}
We consider a coupled system of $q$-fractional differential equations involving two distinct fractional orders and nonlinear $q$-Mittag--Leffler source terms. Throughout this section, let $i=1,2$ and assume that

\begin{equation}
\begin{aligned}
0<\beta_i<\alpha_i<1,&\qquad
\mu_i:=\alpha_i-\beta_i\in(0,1),\\
\theta_i>0,\quad \lambda_i\in\mathbb R,&\qquad
0\le\gamma_i<\mu_i.
\end{aligned}
\label{eq:parameters}
\end{equation}

The restriction $\gamma_i<\mu_i$ is particularly relevant for the well-posedness of the weakly singular contribution arising from the lower-order fractional derivative. We consider the following coupled initial-value problem on $t\in(0,1]$:

\begin{subequations}\label{eq:system}
\begin{equation}\label{eq:system1}
\begin{aligned}
\CD{\alpha_1}\omega_1(t)
+\theta_1t^{-\gamma_1}\CD{\beta_1}\omega_1(t)
=f_1\bigl(t,\omega_1(t),\omega_2(t)\bigr)
+\lambda_1 E_{q,\alpha_1}\bigl(\omega_1(t)\bigr).
\end{aligned}
\end{equation}
\begin{equation}\label{eq:system2}
\begin{aligned}
\CD{\alpha_2}\omega_2(t)
+\theta_2t^{-\gamma_2}\CD{\beta_2}\omega_2(t)
=f_2\bigl(t,\omega_1(t),\omega_2(t)\bigr)
+\lambda_2 E_{q,\alpha_2}\bigl(\omega_2(t)\bigr).
\end{aligned}
\end{equation}
\end{subequations}

subject to the initial conditions
\begin{equation}
\omega_1(0)=a_1,
\qquad
\omega_2(0)=a_2.
\label{eq:initial}
\end{equation}
Here, $f_i:[0,1]\times\mathbb{R}^2\to\mathbb{R}$ are continuous nonlinear coupling functions, while $E_{q,\alpha_i}$ denotes the corresponding $q$-Mittag--Leffler function. The latter provides a genuinely non-polynomial source mechanism and allows the model to incorporate the characteristic response of $q$-fractional systems beyond the usual polynomial or elementary nonlinearities. In particular, $\lambda_i=0$ recovers the associated coupled system without the $q$-Mittag--Leffler contribution. The structure of \eqref{eq:system} contains two distinct sources of nonlocality. The fractional derivatives introduce memory effects, whereas the factors $t^{-\gamma_i}$ generate a weak singularity at the initial point at this problem.
Let us first formulate the given problem into an equivalent Volterra integral system. For this purpose, introduce
\begin{equation*}
g_i(t):=\CD{\alpha_i}\omega_i(t),
\qquad i=1,2,
\end{equation*}
and suppose that $g_i\in C([0,1])$. Since $[0,1]$ is compact, each
$g_i$ is bounded. Using the fractional-order decomposition
$\alpha_i=\beta_i+\mu_i$ and the composition property
\eqref{eq:lowercomposition}, we obtain
\begin{equation*}
\CD{\beta_i}\omega_i(t)
=
\Jq{\mu_i}g_i(t).
\label{eq:lowerfromhighest}
\end{equation*}
Accordingly, define the nonlinear operator

\begin{equation}\label{eq:Gdef}
\begin{aligned}
\mathcal G_i[\omega_1,\omega_2](t)
:={}&-\theta_i t^{-\gamma_i}\CD{\beta_i}\omega_i(t)+f_i\bigl(t,\omega_1(t),\omega_2(t)\bigr)+\lambda_iE_{q,\alpha_i}\bigl(\omega_i(t)\bigr).
\end{aligned}
\end{equation}

for $i=1,2$. The original system can then be viewed as
\[
\CD{\alpha_i}\omega_i
=
\mathcal{G}_i[\omega_1,\omega_2],
\qquad i=1,2.
\]
This representation isolates the nonlinear coupling, the weakly singular fractional interaction, and the $q$-Mittag--Leffler source within a single operator and forms the basis for the subsequent fixed-point and numerical analysis.
\begin{lemma}[{\bf Equivalent Volterra-type Formulation}]
\label{lem:volterra}
Assume that the operators involved in Lemma~\ref{lem:composition} exist. Then the initial-value problem
\eqref{eq:system}--\eqref{eq:initial} is equivalent, in the corresponding operator domain, to the Volterra-type integral system
\begin{equation}
\omega_i(t)
=
a_i+\Jq{\alpha_i}
\mathcal{G}_i[\omega_1,\omega_2](t),
\qquad i=1,2.
\label{eq:volterra}
\end{equation}
Moreover, the lower-order fractional derivatives are represented by
\begin{equation}
\CD{\beta_i}\omega_i(t)
=
\Jq{\mu_i}
\mathcal{G}_i[\omega_1,\omega_2](t),
\qquad i=1,2.
\label{eq:lowerVolterra}
\end{equation}
\end{lemma}

\begin{proof}
From \eqref{eq:system} and \eqref{eq:Gdef}, we have
\[
\CD{\alpha_i}\omega_i(t)
=
\mathcal{G}_i[\omega_1,\omega_2](t).
\]
Applying $\Jq{\alpha_i}$ to both sides and invoking the composition property \eqref{eq:composition} gives
\[
\omega_i(t)-\omega_i(0)
=
\Jq{\alpha_i}
\mathcal{G}_i[\omega_1,\omega_2](t).
\]
Using the initial condition \eqref{eq:initial} yields \eqref{eq:volterra}. Similarly, applying $\Jq{\mu_i}$ and using \eqref{eq:lowercomposition} gives \eqref{eq:lowerVolterra}. Conversely, evaluating \eqref{eq:volterra} at $t=0$ yields $\omega_i(0)=a_i$. Applying $\CD{\alpha_i}$ to \eqref{eq:volterra} and using \eqref{eq:composition} recovers
\[
\CD{\alpha_i}\omega_i(t)
=
\mathcal{G}_i[\omega_1,\omega_2](t),
\]
which, together with \eqref{eq:Gdef}, is precisely \eqref{eq:system}. Hence the two formulations are equivalent.
\end{proof}

\begin{remark}
The $q$-integral representation \eqref{eq:volterra} provides a convenient framework for establishing the existence and uniqueness of solutions to the given problem \eqref{eq:system}--\eqref{eq:initial}. In particular, by defining an appropriate operator associated with the right-hand side of \eqref{eq:volterra}, the problem can be reformulated as a fixed-point problem in a suitable function space. The existence of a solution can then be established by applying an appropriate fixed-point theorem, while uniqueness follows from the corresponding contraction principle. 
These representations serve as the basis for the subsequent existence and uniqueness analysis.
\end{remark}

\section{Quantum Calculus adopted Wavelet Discretization}\label{Section_Methodology}
The present section develops a Haar wavelet collocation framework adapted to the underlying \(q\)-calculus setting for the numerical approximation of the considered coupled fractional differential system. Unlike conventional Haar wavelet schemes formulated in the classical calculus framework, the proposed approach incorporates the relevant quantum-calculus operators into the discretization procedure, thereby preserving the mathematical structure of the governing \(q\)-fractional model. The detailed formulation of the proposed method, including the Haar wavelet basis, the operational relations required for discretization, and the resulting collocation system, is presented below. Let
\begin{equation}
	g_{1,L}(t)=\sum_{j=1}^{N_L}c_jh_j(t),
	\qquad
	g_{2,L}(t)=\sum_{j=1}^{N_L}d_jh_j(t),
	\label{eq:gexpansion}
\end{equation}
where $c=(c_1,\dots,c_{N_L})^T$ and $d=(d_1,\dots,d_{N_L})^T$ are unknown coefficient vectors.
	
\subsection{\texorpdfstring{$q$}{q}-fractional Reconstruction}
For $\nu>0$, the $q$-fractional integral of the Haar wavelet
$h_j$ is defined by
\begin{equation}
\begin{aligned}
\rho_{j,\nu}^{(q)}(t)
:=\Jq{\nu}h_j(t)=\frac{1}{\Gamma_q(\nu)}
\int_0^t(t-qs)_q^{\nu-1}h_j(s)\dq s.
\end{aligned}
\label{eq:rho}
\end{equation}
Accordingly, the reconstructed approximations to the solution components and their lower-order fractional derivatives are given by

\begin{subequations}\label{eq:reconstruction}
\begin{equation}\label{eq:reconstruction1}
\begin{aligned}
\omega_{1,L}(t)
&=a_1+\sum_{j=1}^{N_L}c_j\rho_{j,\alpha_1}^{(q)}(t),\\
v_{1,L}(t)
&:=\sum_{j=1}^{N_L}c_j\rho_{j,\mu_1}^{(q)}(t).
\end{aligned}
\end{equation}
\begin{equation}\label{eq:reconstruction2}
\begin{aligned}
\omega_{2,L}(t)
&=a_2+\sum_{j=1}^{N_L}d_j\rho_{j,\alpha_2}^{(q)}(t),\\
v_{2,L}(t)
&:=\sum_{j=1}^{N_L}d_j\rho_{j,\mu_2}^{(q)}(t).
\end{aligned}
\end{equation}
\end{subequations}

Here, $c_j$ and $d_j$ are unknown expansion coefficients, and $v_{i,L}$ approximates $\CD{\beta_i}\omega_i$. No additive constant is required in $v_{i,L}$, since the lower-order fractional derivatives are determined by the composition relation \eqref{eq:lowercomposition}. Substituting the defining Jackson sum into \eqref{eq:rho} gives the identity
\begin{equation}
		\rho_{j,\nu}^{(q)}(t)
		=\frac{(1-q)t}{\Gamma_q(\nu)}
		\sum_{n=0}^{\infty}q^n
		\bigl(t-q^{n+1}t\bigr)_q^{\nu-1}
		h_j(tq^n).
		\label{eq:rhoseries}
\end{equation}
In particular, for the scaling function $h_1(t)\equiv 1$, the above expression simplifies to $\rho_{1,\nu}^{(q)}(t)=\frac{t^\nu}{\Gamma_q(\nu+1)}$. Define the collocation points $t_r=\left(r-\frac12\right)h_L,~~ r=1,\dots,N_L.$ These are strictly positive for each fixed $L$, so no collocation equation is imposed at the singular endpoint $t=0$. Substituting \eqref{eq:gexpansion}--\eqref{eq:reconstruction} into \eqref{eq:system} gives the residuals:

\begin{subequations}\label{eq:residuals}
\begin{equation}\label{eq:R1}
\begin{aligned}
R_{1,r}(c,d)\coloneqq{}&
\sum_{j=1}^{N_L}c_j h_j(t_r)+\theta_1t_r^{-\gamma_1}v_{1,L}(t_r)-f_1\bigl(t_r,\omega_{1,L}(t_r),\qquad\omega_{2,L}(t_r)\bigr)-\lambda_1 E_{q,\alpha_1}\bigl(\omega_{1,L}(t_r)\bigr).
\end{aligned}
\end{equation}
\begin{equation}\label{eq:R2}
\begin{aligned}
R_{2,r}(c,d)\coloneqq{}&
\sum_{j=1}^{N_L}d_j h_j(t_r)+\theta_2t_r^{-\gamma_2}v_{2,L}(t_r)-f_2\bigl(t_r,\omega_{1,L}(t_r),\qquad\omega_{2,L}(t_r)\bigr)-\lambda_2 E_{q,\alpha_2}\bigl(\omega_{2,L}(t_r)\bigr).
\end{aligned}
\end{equation}
\end{subequations}

Therefore, the exact collocation system reduces to
\begin{equation}
	R_{i,r}(c,d)=0,
	\qquad i=1,2,
	\quad r=1,\dots,N_L.
	\label{eq:collocationsystem}
\end{equation}

\section{Qualitative Analysis}\label{Section_Qualitative}
This section investigates the existence and uniqueness of solutions to the proposed coupled $q$-fractional system. By exploiting the equivalent Volterra-type integral formulation established in Lemma~\ref{lem:volterra}, we derive sufficient conditions for the well-posedness of the initial-value problem. The analysis relies on an appropriate fixed-point theorem and suitable assumptions on the nonlinear operators involved.	
\begin{lemma}\label{lem:Jbound}
	Let $\nu>0$ and let $v:[0,1]\to\R$ be bounded. Then $\abs{\Jq{\nu}v(t)} \le C_\nu t^\nu\norm{v}_\infty, 0\le t\le1$, and therefore 
    $$\norm{\Jq{\nu}v}_\infty\le C_\nu\norm{v}_\infty,$$ 
    where $C_\nu:=\frac{1}{\Gamma_q(\nu+1)}$.
\end{lemma}
	
\begin{proof}
		The kernel in \eqref{eq:qfracint} is nonnegative for $0\le s\le t$. Hence
		$
			\abs{\Jq{\nu}v(t)}
			\le\frac{\norm{v}_\infty}{\Gamma_q(\nu)}
			\int_0^t(t-qs)_q^{\nu-1}\dq s.
	$
		The $q$-beta identity gives
		$
			\int_0^t(t-qs)_q^{\nu-1}\dq s
			=\frac{\Gamma_q(\nu)}{\Gamma_q(\nu+1)}t^\nu.
		$
		Substitution proves $\abs{\Jq{\nu}v(t)}
		\le C_\nu t^\nu\norm{v}_\infty$, and $\norm{\Jq{\nu}v}_\infty\le C_\nu\norm{v}_\infty$ follows because $t^\nu\le1$.
\end{proof}
	
\begin{lemma}\label{lem:Jcontinuity}
Let $v\in \mathcal{C}([0,1])$ and $\nu>0$. Then $J_q^\nu v\in \mathcal{C}([0,1]).$
\end{lemma}

\begin{proof}
Based on the Jackson-series representation of the $q$-fractional integral (see Definition~\ref{Def_q_Integral} and and \cite{KacCheung2002,Ernst2012}) with its scaling property, one obtain
\[J_q^\nu v(t)=\sum_{n=0}^{\infty}S_n(t),\qquad 0\leq t\leq1.\]
Since $v\in \mathcal{C}([0,1])$, it is bounded and hence $|v(tq^n)|\leq \|v\|_\infty$. Consequently, $|S_n(t)|\leq\|v\|_\infty M_n,$ where $M_n$ is the corresponding nonnegative term obtained from the Jackson representation after setting $t=1$ and replacing $v$ by the constant function $1$. By Lemma~\ref{lem:Jbound},
$
\sum_{n=0}^{\infty}M_n
=
J_q^\nu 1(1)
<\infty.
$
Therefore the series defining $J_q^\nu v$ converges uniformly on
$[0,1]$ by the Weierstrass $M$-test. Each term is continuous for
$t>0$, and hence the sum is continuous for $t>0$. Finally, Lemma~\ref{lem:Jbound}
gives
\[
|J_q^\nu v(t)|
\leq
C_\nu t^\nu\|v\|_\infty
\longrightarrow0
\qquad (t\to0^+),
\]
which proves continuity at $t=0$.
\end{proof}

	\begin{lemma}
		\label{lem:Wbound}
		Let $\mu>0$ and $0\le\gamma<\mu$. Define
		\[
			(\mathcal W_{\mu,\gamma}v)(t)=
			\begin{cases}
				t^{-\gamma}\Jq{\mu}v(t), & t>0,\\
				0, & t=0.
			\end{cases}
			\label{eq:Wdef}
		\]
		Then $\mathcal W_{\mu,\gamma}$ is linear and bounded on the space of bounded functions, with
		$
			\abs{(\mathcal W_{\mu,\gamma}v)(t)}
			\le C_\mu t^{\mu-\gamma}\norm{v}_\infty$ and $
			\norm{\mathcal W_{\mu,\gamma}v}_\infty
			\le C_\mu\norm{v}_\infty.
			\label{eq:Wbound}
		$ If $v$ is continuous, then $\mathcal W_{\mu,\gamma}v$ is continuous at $t=0$.
	\end{lemma}
	
\begin{proof}
The operator $W_{\mu,\gamma}$ is linear because the
$q$-fractional integral $J_q^\mu$ is linear. Indeed, for
$a,b\in\mathbb{R}$ and bounded functions $u,v$, we have, for $t>0$,

\[
\begin{aligned}
W_{\mu,\gamma}(au+bv)(t)
&=t^{-\gamma}J_q^\mu(au+bv)(t)=a\,t^{-\gamma}J_q^\mu u(t)+b\,t^{-\gamma}J_q^\mu v(t)=aW_{\mu,\gamma}u(t)+bW_{\mu,\gamma}v(t).
\end{aligned}
\]

At $t=0$, both sides are equal to zero by definition. Hence
$W_{\mu,\gamma}$ is linear. For $t>0$, multiplying the estimate in Lemma~\ref{lem:Jbound} by $t^{-\gamma}$
gives
$
|W_{\mu,\gamma}v(t)|
=
t^{-\gamma}|J_q^\mu v(t)|
\leq
C_\mu t^{\mu-\gamma}\|v\|_\infty.
$
Since $\mu-\gamma>0$ and $0\leq t\leq1$, it follows that
$
|W_{\mu,\gamma}v(t)|
\leq
C_\mu\|v\|_\infty,
$
and therefore
$
\|W_{\mu,\gamma}v\|_\infty
\leq
C_\mu\|v\|_\infty.
$
Thus, $W_{\mu,\gamma}$ is bounded. Further, if $v$ is continuous, then $J_q^\mu v$ is continuous by
Lemma~\ref{lem:Jcontinuity}. Moreover,
\[
|W_{\mu,\gamma}v(t)|
\leq
C_\mu t^{\mu-\gamma}\|v\|_\infty
\longrightarrow 0
\qquad\text{as } t\to0^+.
\]
Since $W_{\mu,\gamma}v(0)=0$, the function
$W_{\mu,\gamma}v$ is also continuous at $t=0$.
\end{proof}
	
To establish the existence and uniqueness of a solution, we formulate the equivalent Volterra-type system as a fixed-point problem in a suitable Banach space. The corresponding nonlinear operator is constructed in terms of the highest-order fractional derivatives, with the solution components recovered through the \(q\)-fractional integral operators. Let
	$
		\mathcal B=\{v:[0,1]\to\R:\norm{v}_\infty<\infty\}$ for 
		$\mathcal X=\mathcal B\times\mathcal B$ with $\norm{g}_{\mathcal X}=\max\{\norm{g_1}_\infty,\norm{g_2}_\infty\}.
$
	For $R>0$, set
	$
		\mathbb B_R=\{g\in\mathcal X:\norm{g}_{\mathcal X}\le R\}.
	$
	For $g=(g_1,g_2)\in\mathcal X$, define
	$
		\Omega_i[g](t)=a_i+\Jq{\alpha_i}g_i(t),
		i=1,2,
		\label{eq:Omega}
	$
	and define $\mathcal F=(\mathcal F_1,\mathcal F_2)$ by

\[
\begin{aligned}
\mathcal F_i(g)(t)={}&-\theta_i\mathcal W_{\mu_i,\gamma_i}g_i(t)+f_i\bigl(t,\Omega_1[g](t),\qquad\Omega_2[g](t)\bigr)+\lambda_iE_{q,\alpha_i}\bigl(\Omega_i[g](t)\bigr),
\quad i=1,2.
\end{aligned}
\]

A fixed point $g=\mathcal F(g)$ is precisely the system written in terms of the highest-order derivatives $g_i=\CD{\alpha_i}\omega_i$. For fixed $R>0$ and numbers $r_1,r_2>0$, define

\[
\begin{aligned}
\mathcal Q_R &=[0,1]\times[-r_1,r_1]\times[-r_2,r_2],\\
M_{f_i}(R)&:=\sup_{(t,x_1,x_2)\in\mathcal Q_R}
\lvert f_i(t,x_1,x_2)\rvert,\hspace{5em}i=1,2.
\end{aligned}
\]

Let $L_{f_i}$ be a Lipschitz constant of $f_i$ with respect
to its state variables on $\mathcal{Q}_R$, that is,
\[
\begin{aligned}
&|f_i(t,x_1,x_2)-f_i(t,y_1,y_2)|\le L_{f_i}(|x_1-y_1|+|x_2-y_2|).
\end{aligned}
\]
Furthermore, define

\[
\begin{aligned}
M_{E_i}(R)
&:=\sup_{|x|\le r_i}|E_{q,\alpha_i}(x)|,\\
L_{E_i}(R)
&:=\sup_{|x|\le r_i}
\left|E'_{q,\alpha_i}(x)\right|.
\end{aligned}
\]

\begin{assumption}[Sufficient conditions for contractivity]
\label{ass:contraction}
There exist constants $R>0$ and $r_i>0$, $i=1,2$, such that
\begin{equation}
    |a_i|+C_{\alpha_i}R
    \le r_i
    <(1-q)^{-\alpha_i},
    \qquad i=1,2.
    \label{eq:statebound}
\end{equation}
Moreover, the nonlinear operator satisfies the invariance conditions
\begin{equation}
    \theta_i C_{\mu_i}R
    +M_{f_i}(R)
    +|\lambda_i|M_{E_i}(R)
    \le R,
    \qquad i=1,2,
    \label{eq:invariant}
\end{equation}
and the contraction condition
\begin{equation*}
    \kappa:=\max\{\kappa_1,\kappa_2\}<1,
\end{equation*}
where

\begin{equation*}
\begin{aligned}
\kappa_i:={}&\theta_i C_{\mu_i}+L_{f_i}(C_{\alpha_1}+C_{\alpha_2})+|\lambda_i|L_{E_i}(R)C_{\alpha_i},
\quad i=1,2.
\end{aligned}
\end{equation*}

Here, $C_{\alpha_i}$ and $C_{\mu_i}$ denote the corresponding
$q$-fractional integral bounds, while $M_{f_i}(R)$,
$M_{E_i}(R)$, $L_{f_i}$, and $L_{E_i}(R)$ are the bounds and
Lipschitz constants specified above.
\end{assumption}
These conditions are sufficient  to ensure that the nonlinear operator maps the closed ball into itself which is also a contraction. They are not asserted to be necessary conditions. The following theorem establishes the existence and uniqueness of the solution by applying the Banach fixed-point theorem to the operator \(\mathcal F\). 
	
\begin{theorem}[{\bf Existence and Uniqueness}]
		\label{thm:existence}
		Under \eqref{eq:parameters} and Assumption~\ref{ass:contraction}, the operator $\mathcal F$ maps $\mathbb B_R$ into itself and is a contraction with constant $\kappa$. Consequently, there exists a unique fixed point $g=(g_1,g_2)\in\mathbb B_R$. This fixed point is continuous. The reconstructed functions
		\begin{equation*}
			\omega_i(t)=a_i+\Jq{\alpha_i}g_i(t),
			\qquad i=1,2,
		\end{equation*}
form the unique solution of the integral representation
\eqref{eq:volterra} in this bounded class. Whenever the composition identities of Lemma~\ref{lem:composition} apply, they also form the unique solution of \eqref{eq:system}--\eqref{eq:initial} in the same class.
\end{theorem}
	
\begin{proof}
Let $g\in\mathbb B_R$. Equation (\ref{eq:statebound}) and Lemma~\ref{lem:Jbound} give
		$
			\abs{\Omega_i[g](t)}
			\le\abs{a_i}+C_{\alpha_i}\norm{g_i}_\infty
			\le\abs{a_i}+C_{\alpha_i}R
			\le r_i.
	$
		Therefore every state argument remains in the compact sets used in Assumption~\ref{ass:contraction}. Lemma~\ref{lem:Wbound} and \eqref{eq:invariant} yield
		$
			\norm{\mathcal F_i(g)}_\infty
			\le\theta_iC_{\mu_i}R+M_{f_i}(R)+\abs{\lambda_i}M_{E_i}(R)
			\le R.
		$
		Thus $\mathcal F(\mathbb B_R)\subseteq\mathbb B_R$. Lemmas~\ref{lem:Jcontinuity} and \ref{lem:Wbound}, continuity of $f_i$, and Proposition~\ref{prop:qMLdisk} also show that $\mathcal F$ maps $\mathcal{C}([0,1])^2\cap\mathbb B_R$ into itself. Let $g,v\in\mathbb B_R$. Then, by Lemmas~\ref{lem:Jbound} and \ref{lem:Wbound}, we yield
\begin{align}
\norm{\mathcal F_i(g)-\mathcal F_i(v)}_\infty
\le{}&\theta_iC_{\mu_i}\norm{g_i-v_i}_\infty
+L_{f_i}\bigl(C_{\alpha_1}\norm{g_1-v_1}_\infty
+C_{\alpha_2}\norm{g_2-v_2}_\infty\bigr)
+|\lambda_i|L_{E_i}(R)C_{\alpha_i}
\nonumber\\
&\qquad\times\norm{g_i-v_i}_\infty
\nonumber\\
\le{}&\kappa_i\norm{g-v}_{\mathcal X}.
\nonumber
\end{align}

		Taking the maximum over $i$ gives
		$
			\norm{\mathcal F(g)-\mathcal F(v)}_{\mathcal X}
			\le\kappa\norm{g-v}_{\mathcal X}.
		$
		The Banach fixed-point theorem yields a unique fixed point in $\mathbb B_R$. The closed set $C([0,1])^2\cap\mathbb B_R$ is complete, invariant, and subject to the same contraction; it therefore has a fixed point, which by uniqueness in $\mathbb B_R$ is the same fixed point. Hence $g$ is continuous. The reconstruction and Lemma~\ref{lem:composition} give the stated integral and differential solutions.
	\end{proof}

\section{Stability and Convergence Analysis}\label{Section_Convergence}
Following the existence and uniqueness results, this section examines the stability and convergence properties of the proposed \(q\)-Haar wavelet collocation method. We investigate the sensitivity of the numerical solution to perturbations and assess its convergence toward the solution of the underlying \(q\)-fractional system. The analysis provides a theoretical basis for the reliability of the proposed numerical scheme and complements the numerical results presented in the subsequent section. The following theorem establishes the stability of the fixed-point solution with respect to both residual errors and perturbations in the underlying nonlinear operator.	
\begin{theorem}\label{thm:stability}
Let $g\in\mathbb B_R$ be the unique fixed point of $\mathcal F$ guaranteed by Theorem~\ref{thm:existence}. Then the following stability estimates hold:
\begin{enumerate}[label=\textup{(\roman*)}]
\item For every $\widehat g\in\mathbb B_R$, the error between the exact fixed point and the approximate solution satisfies
	\begin{equation}
		\norm{g-\widehat g}_{\mathcal X}
		\le\frac{1}{1-\kappa}
		\norm{\widehat g-\mathcal F(\widehat g)}_{\mathcal X}.
		\label{eq:residualstability}
	\end{equation}
			
\item Let $\widetilde{\mathcal F}:\mathbb B_R\to\mathbb B_R$ be a contraction, and let $\widetilde g\in\mathbb B_R$ be its fixed point. If
			\begin{equation*}
				\Delta_F:=\sup_{v\in\mathbb B_R}
				\norm{\mathcal F(v)-\widetilde{\mathcal F}(v)}_{\mathcal X},
			\end{equation*}
			then
			\begin{equation}
				\norm{g-\widetilde g}_{\mathcal X}
				\le\frac{\Delta_F}{1-\kappa}.
				\label{eq:operatorstability}
			\end{equation}
		\end{enumerate}
\end{theorem}
	
\begin{proof}
	For part (i), using $g=\mathcal F(g)$, we have
		\begin{align*}
			\norm{g-\widehat g}_{\mathcal X}
			&\le\norm{\mathcal F(\widehat g)-\widehat g}_{\mathcal X}+\norm{\mathcal F(g)-\mathcal F(\widehat g)}_{\mathcal X}
			\\
			&\le\norm{\mathcal F(\widehat g)-\widehat g}_{\mathcal X}+\kappa\norm{g-\widehat g}_{\mathcal X}
			.
		\end{align*}
		Rearranging proves \eqref{eq:residualstability}. To prove part (ii), we use  $\widetilde g=\widetilde{\mathcal F}(\widetilde g)$ and $g=\mathcal F(g)$ :
		\begin{align*}
			\norm{g-\widetilde g}_{\mathcal X}
			&\le\norm{\mathcal F(g)-\mathcal F(\widetilde g)}_{\mathcal X}
			+\norm{\mathcal F(\widetilde g)-\widetilde{\mathcal F}(\widetilde g)}_{\mathcal X}\le\kappa\norm{g-\widetilde g}_{\mathcal X}+\Delta_F.
		\end{align*}
		Rearranging gives \eqref{eq:operatorstability}.
	\end{proof}
	
\begin{remark}
	Estimate \eqref{eq:residualstability} controls the effect of terminating the nonlinear solver with a nonzero residual. Estimate \eqref{eq:operatorstability} controls perturbations caused, for example, by truncating the Jackson series in \eqref{eq:rhoseries}. Thus discretization, Jackson-sum, and nonlinear-solver errors can be separated rather than conflated.
\end{remark}
Define the midpoint sampling operator $\mathcal I_L:\mathcal B\to V_L$ by $(\mathcal I_Lv)(t)=v(t_r), t\in I_r, r=1,\dots,N_L$. For vectors, $\mathcal I_L$ is applied component-wise. The following lemma establishes its stability and approximation properties in terms of the H\"{o}lder regularity of the underlying function. In particular, it yields a first-order error estimate for functions in $W^{1,\infty}(0,1)$.
\begin{lemma}\label{lem:interpolation}
	For every bounded $v$,
		\begin{equation}
			\norm{\mathcal I_Lv}_\infty\le\norm{v}_\infty.
			\label{eq:Istability}
		\end{equation}
		Let $0<\sigma\le1$, and suppose $v\in C^{0,\sigma}([0,1])$ i.e., $v$ is continuous with a bounded norm 
		\begin{equation*}
			[v]_{\sigma}:=\sup_{s\ne t}\frac{\abs{v(t)-v(s)}}{\abs{t-s}^{\sigma}}.
		\end{equation*}
		Then
		\begin{equation}
			\norm{v-\mathcal I_Lv}_\infty
			\le [v]_{\sigma}\left(\frac{\tau_L}{2}\right)^{\sigma}.
			\label{eq:Ierror}
		\end{equation}
		In particular, if $v\in W^{1,\infty}(0,1)$, then \eqref{eq:Ierror} holds with $\sigma=1$ and $[v]_1\le\norm{v'}_\infty$.
	\end{lemma}
	
	\begin{proof}
		The first assertion follows because $\mathcal I_Lv$ assumes only sampled values of $v$. If $t\in I_r$, then
		\begin{equation*}
			\abs{v(t)-v(t_r)}
			\le [v]_{\sigma}\abs{t-t_r}^{\sigma}
			\le [v]_{\sigma}\left(\frac{\tau_L}{2}\right)^{\sigma}.
		\end{equation*}
		Taking the supremum proves \eqref{eq:Ierror}. The last assertion follows from the Lipschitz representative of a $W^{1,\infty}$ function.
	\end{proof}
	\noindent
	
	\begin{proposition}
		\label{prop:equivalence}
		Define the exact discrete operator
	$
		\mathcal T_L=\mathcal I_L\mathcal F
		\label{eq:TL}
	$. Let $g_L=(g_{1,L},g_{2,L})\in V_L\times V_L$ and reconstruct $\omega_{i,L}$ and $v_{i,L}$ by \eqref{eq:reconstruction}, using exact $q$-fractional Haar integrals. Then
		$
			g_L=\mathcal T_Lg_L
		$
		if and only if the midpoint residual equations \eqref{eq:collocationsystem} hold.
	\end{proposition}
	
	\begin{proof}
Both $g_{i,L}$ and $\mathcal{I}_L\mathcal{F}_i(g_L)$ belong to
$V_L$ and are therefore constant on every dyadic cell $I_r$.
Consequently,
$
g_{i,L}
=
\mathcal{I}_L\mathcal{F}_i(g_L)
$
if and only if their values coincide at the midpoint $t_r$ of each
cell, that is,
$
g_{i,L}(t_r)
=
\mathcal{F}_i(g_L)(t_r),
r=1,\ldots,N_L.
$
After substituting the Haar representations \eqref{eq:gexpansion} and the
reconstruction formulas (\ref{eq:reconstruction1})-(\ref{eq:reconstruction2}), these midpoint identities are exactly
the residual equations (\ref{eq:R1}), (\ref{eq:R2}), and (\ref{eq:collocationsystem}).
\end{proof}

\begin{remark}
The convergence estimate below uses the approximation property of the midpoint projection onto the piecewise-constant space $V_L$. The classical Haar basis is used to represent elements of $V_L$ and to construct the algebraic collocation scheme; however, no additional
orthogonality property of the individual Haar functions is required in the proof of the following error estimate.
\end{remark}

\subsection{Error Estimates}
\begin{theorem}\label{thm:convergence}
		Assume \eqref{eq:parameters} and Assumption~\ref{ass:contraction}. Then $\mathcal T_L$ has a unique fixed point
		$
			g_L\in\mathbb B_R\cap(V_L\times V_L),
		$
		which is the unique exact midpoint collocation solution. Let $g$ be the continuous fixed point from Theorem~\ref{thm:existence}. If, for some $0<\sigma\le1$,
		\begin{equation*}
			g_i\in C^{0,\sigma}([0,1]),
			\qquad
			H_{\sigma}:=\max_{i=1,2}[g_i]_{\sigma},
		\end{equation*}
		then
		\begin{equation}
			\norm{g-g_L}_{\mathcal X}
			\le\frac{H_{\sigma}}{1-\kappa}
			\left(\frac{\tau_L}{2}\right)^{\sigma}.
			\label{eq:gconvergence}
		\end{equation}
		Moreover, with
		\begin{equation*}
			\omega_{i,L}=a_i+\Jq{\alpha_i}g_{i,L},
			\qquad
			v_{i,L}=\Jq{\mu_i}g_{i,L},
		\end{equation*}
		we have
		\begin{align}
			\norm{\omega_i-\omega_{i,L}}_\infty
			&\le\frac{C_{\alpha_i}H_{\sigma}}{1-\kappa}
			\left(\frac{\tau_L}{2}\right)^{\sigma},
			\label{eq:omegaconvergence}\\
			\norm{\CD{\beta_i}\omega_i-v_{i,L}}_\infty
			&\le\frac{C_{\mu_i}H_{\sigma}}{1-\kappa}
			\left(\frac{\tau_L}{2}\right)^{\sigma},
			\label{eq:vconvergence}\\
			\norm{t^{-\gamma_i}\bigl(\CD{\beta_i}\omega_i-v_{i,L}\bigr)}_\infty
			&\le\frac{C_{\mu_i}H_{\sigma}}{1-\kappa}
			\left(\frac{\tau_L}{2}\right)^{\sigma},
			\label{eq:singularconvergence}
		\end{align}
		where the expression in \eqref{eq:singularconvergence} is extended by zero at $t=0$. Consequently,
		\begin{equation*}
			\norm{\omega_i-\omega_{i,L}}_\infty
			=\mathcal{O}(\tau_L^{\sigma})=\mathcal{O}\!\left(2^{-\sigma(L+1)}\right).
		\end{equation*}
		If $g_i\in W^{1,\infty}(0,1)$, then $\sigma=1$ is admissible and the theoretical rate is first order:
		\begin{equation*}
			\norm{\omega_i-\omega_{i,L}}_\infty=\mathcal{O}(\tau_L)=\mathcal{O}(2^{-(L+1)}).
		\end{equation*}
	\end{theorem}
	
	\begin{proof}
		The role of the Haar discretization enters through the piecewise-constant space $V_L$ and the midpoint projection $\mathcal{I}_L$. The error estimate itself follows from the
approximation property of $\mathcal{I}_L$ together with the contraction property of $\mathcal{F}$. By \eqref{eq:Istability} and the self-map property of $\mathcal F$, the operator $\mathcal T_L=\mathcal I_L\mathcal F$ maps $\mathbb B_R$ into itself. Moreover,
		$
			\norm{\mathcal T_Lg-\mathcal T_Lv}_{\mathcal X}
			\le\norm{\mathcal F(g)-\mathcal F(v)}_{\mathcal X}
			\le\kappa\norm{g-v}_{\mathcal X}.
		$
		Hence $\mathcal T_L$ is a contraction and has a unique fixed point in $\mathbb B_R$. Since its range is contained in $V_L\times V_L$, the fixed point belongs to that space. Proposition~\ref{prop:equivalence} identifies it with the unique exact collocation solution. By Assumption~\ref{ass:contraction} and Theorem~\ref{thm:existence},
$\mathcal{F}$ maps $\mathbb{B}_R$ into itself and satisfies
\[
\|\mathcal{F}(u)-\mathcal{F}(v)\|_{\mathcal X}
\leq
\kappa\|u-v\|_{\mathcal X},
\qquad
0\leq\kappa<1.
\]
Hence the factor $1-\kappa$ appearing below is strictly positive. Because $g=\mathcal F(g)$ and $g_L=\mathcal I_L\mathcal F(g_L)$,
		\begin{align*}
\norm{g-g_L}_{\mathcal X}
&\le\norm{g-\mathcal I_Lg}_{\mathcal X}\\
&\quad+\norm{\mathcal I_L\mathcal F(g)-\mathcal I_L\mathcal F(g_L)}_{\mathcal X}\\
&\le\norm{g-\mathcal I_Lg}_{\mathcal X}
+\kappa\norm{g-g_L}_{\mathcal X}.
\end{align*}
		Therefore,
		$
			(1-\kappa)\norm{g-g_L}_{\mathcal X}
			\le\norm{g-\mathcal I_Lg}_{\mathcal X}.
		$
		Applying \eqref{eq:Ierror} componentwise gives \eqref{eq:gconvergence}.\\
		
		\noindent
		The reconstruction formulas give
		$
			\omega_i-\omega_{i,L}=\Jq{\alpha_i}(g_i-g_{i,L})$ and 
			$
			\CD{\beta_i}\omega_i-v_{i,L}=\Jq{\mu_i}(g_i-g_{i,L}).
		$
		Hence, equation \eqref{eq:gconvergence} and Lemmas~\ref{lem:Jbound},~\ref{lem:Wbound}, prove \eqref{eq:omegaconvergence}--\eqref{eq:singularconvergence}. Since $\tau_L=2^{-(L+1)}$, the order statement follows.
	\end{proof}

\begin{remark}\label{remark_Convergence}
	Theorem~\ref{thm:convergence} gives convergence and stability of the midpoint collocation scheme under the contraction assumption. By \eqref{eq:gconvergence}, the rate is governed by the H\"{o}lder regularity of the fixed point $g$: if $g_i\in C^{0,\sigma}([0,1])$, the error is $\mathcal{O}(\tau_L^\sigma)$ (cf.\ \cite{Majak2015}), with a constant depending on the Hölder seminorm and the contraction factor $\kappa$ and first-order convergence is guaranteed when $g_i\in W^{1,\infty}(0,1)$. However, the numerical results in the subsequent section exhibit first-order convergence even when $g_i\notin W^{1,\infty}(0,1)$. The estimates for $\omega_i$ and $v_{i,L}$ show that this rate carries over to the reconstructed quantities, and \eqref{eq:singularconvergence} keeps the weighted error uniformly controlled despite the singular weight at $t=0$.
\end{remark}

\section{Damped Newton Iteration}\label{Sec_Newton}
Set $X=(c_1,\dots,c_{N_L},d_1,\dots,d_{N_L})^T\in\R^{2N_L}$ and order the residual vector in the same manner. At iteration $k$, solve $\mathsf J(X^{(k)})\Delta X^{(k)}=-R(X^{(k)}),$ and use $X^{(k+1)}=X^{(k)}+\eta_k\Delta X^{(k)}, 0<\eta_k\le1.$ A backtracking line search should select $\eta_k$ so that the residual norm decreases and every evaluated $q$-Mittag--Leffler argument remains strictly inside its convergence disk. Let $A_{i,rj}=\rho_{j,\alpha_i}^{(q)}(t_r), B_{i,rj}=\rho_{j,\mu_i}^{(q)}(t_r)$ and $H_{rj}=h_j(t_r).$ Although the linear part of $R_{1,r}$ contains only the coefficients $c_j$, the residual $R_{1,r}$ also depends on $d_j$ through $\omega_{2,L}$ in the nonlinear function $f_1$. Similarly, $R_{2,r}$ depends on $c_j$ through $\omega_{1,L}$ in $f_2$. Hence, we have
\[
\frac{\partial \omega_{1,L}(t_r)}{\partial c_j}
=A_{1,rj},
\qquad
\frac{\partial \omega_{2,L}(t_r)}{\partial d_j}
=A_{2,rj}.
\] 
Assume that $f_i$ is continuously differentiable in its state variables. The Jacobian entries are then estimated as

\begin{align*}
\frac{\partial R_{1,r}}{\partial c_j}
={}&H_{rj}+\theta_1t_r^{-\gamma_1}B_{1,rj}-\partial_2f_1\,A_{1,rj}-\lambda_1E'_{q,\alpha_1}(\omega_{1,L})A_{1,rj},\\
\frac{\partial R_{1,r}}{\partial d_j}
={}&-\partial_3f_1\,A_{2,rj},\\
\frac{\partial R_{2,r}}{\partial c_j}
={}&-\partial_2f_2\,A_{1,rj},\\
\frac{\partial R_{2,r}}{\partial d_j}
={}&H_{rj}+\theta_2t_r^{-\gamma_2}B_{2,rj}-\partial_3f_2\,A_{2,rj}-\lambda_2E'_{q,\alpha_2}(\omega_{2,L})A_{2,rj}.
\end{align*}

The Jacobian has the block form
\[
	\mathbf{J}
	=\begin{pmatrix}
		\frac{\partial R_{1,r}}{\partial c_j}
		&\frac{\partial R_{1,r}}{\partial d_j}\\
		\frac{\partial R_{2,r}}{\partial c_j}
		&\frac{\partial R_{2,r}}{\partial d_j}
	\end{pmatrix}.
\]	
The derivatives of $f_i$ are evaluated at $(t_r,\omega_{1,L}(t_r),\omega_{2,L}(t_r))$, and the unlabelled Mittag--Leffler arguments are evaluated at the corresponding state value at $t_r$. A numerical solve should be accepted only when both a residual criterion and a step criterion are satisfied. The complete computational procedure is summarized in Algorithm \ref{alg:qHWCM}.
{\small
\begin{figure}[t]
\centering
\begin{minipage}{0.94\textwidth}
\refstepcounter{algorithm}
\noindent\textbf{Algorithm~\thealgorithm. $q$-Haar wavelet collocation method for the coupled $q$-fractional system}\label{alg:qHWCM}\par\medskip\hrule\smallskip
		\begin{algorithmic}[1]
			
			\REQUIRE Resolution level $L$; parameters
			$q,\alpha_i,\beta_i,\theta_i,\gamma_i,\lambda_i$;
			initial values $a_i$; functions $f_i$;
			tolerances $\mathcal{T}_q,\epsilon_R,\epsilon_X$;
			maximum iteration count $k_{\mathrm{max}}$
			
			\ENSURE Approximate solutions $\omega_{1,L}$ and $\omega_{2,L}$,
			or a failure flag
			
			\STATE Set $K = 2^{L}$ and $N = 2K$
			
			\STATE Generate $\tau_{r} = (r - \frac{1}{2}) / N_L$ ,
			$r = 1,\ldots ,N_L$
			
			\STATE Compute $H_{rj} = h_j(\tau_r)$
			
			\STATE Compute the matrices
			$A_{i,rj} = \rho_{j,\alpha_i}^{(q)}(\tau_r)$
			and
			$B_{i,rj} = \rho_{j,\alpha_i - \beta_i}^{(q)}(\tau_r)$
			with Jackson-sum tolerance $\mathcal{T}_q$
			
			\STATE Set $\mathbf{X}^{(0)} = \mathbf{0}$ and $k = 0$
			
			\REPEAT
			
			\STATE Evaluate $\mathbf{R}(\mathbf{X}^{(k)})$
			
			\STATE Assemble $\mathbf{J}(\mathbf{X}^{(k)})$
			
			\STATE Solve $\mathbf{J}\Delta \mathbf{X}^{(k)} = -\mathbf{R}$
			
			\STATE Select $\eta_{k}\in (0,1]$ by backtracking
			
			\STATE Set $\mathbf{X}^{(k + 1)} = \mathbf{X}^{(k)} + \eta_{k}\Delta \mathbf{X}^{(k)}$
			
			\STATE $k\gets k + 1$
			
			\UNTIL{both stopping criteria hold or $k\geq k_{\mathrm{max}}$}
			
			\IF{the stopping criteria do not hold}
			
			\STATE return failure flag
			
			\ENDIF
			
			\STATE Extract $\mathbf{c}^{*}$ and $\mathbf{d}^{*}$ from $\mathbf{X}^{(k)}$
			
			\STATE Reconstruct $\omega_{1,L}$ and $\omega_{2,L}$ from \eqref{eq:reconstruction}
			
			\RETURN $\omega_{1,L},\omega_{2,L}$
			
		\end{algorithmic}
\smallskip\hrule
\end{minipage}
\end{figure}
}

\section{Computational Validation and Insights}\label{Section_Numerical_Examples}
We now assess the performance of the proposed q-HWCM through a series of numerical experiments. All simulations were implemented in Mathematica and run on identical hardware to ensure consistent CPU-time comparisons. The graphs are plotted with the help of MATLAB figure tools. Throughout, errors are measured on a fixed validation grid $\mathcal G\subset[0,1]$ consisting of $2049$ uniformly spaced points, independent of the collocation resolution $L$ used to construct the approximate solution; this decouples the error measurement from the discretization itself and avoids superconvergence artifacts at collocation nodes. For each resolution level $L$, with mesh size $h_L$ and $N=2^{L+1}$ collocation points, we define the maximum absolute error
\[
E_{i,L}
=
\max_{t\in\mathcal G}
\left|
\omega_i(t)-\omega_{i,L}(t)
\right|,
\qquad i=1,2,
\]
where $\omega_i$ denotes the (exact or reference) solution component and $\omega_{i,L}$ its q-HWCM approximation at level $L$. The empirical order of convergence is recovered from two consecutive resolutions via
\[
p_{i,L}
=
\frac{\log(E_{i,L}/E_{i,L+1})}
{\log(h_L/h_{L+1})}
=
\log_2\left(
\frac{E_{i,L}}{E_{i,L+1}}
\right).
\]
the second equality following from the dyadic refinement $h_{L+1}=h_L/2$ inherent to the Haar multiresolution structure. Before reporting results, we verify that the structural hypotheses underlying our convergence theory hold for every case considered. Specifically, Assumption~\ref{ass:contraction} is checked explicitly in each example: fixing the ball radius $R=0.1$, we set
\[
r_i=|a_i|+C_{\alpha_i}R,\qquad i=1,2.
\]
and confirm that the resulting state intervals lie strictly inside the convergence disks of the corresponding $q$-Mittag--Leffler functions. We further verify that the invariant-ball quantities
\[
\Phi_i
=
\theta_iC_{\mu_i}R
+
M_{f_i}(R)
+
|\lambda_i|M_{E_i}(R)
\]
satisfy $\Phi_i<R$, and the contraction constant $\kappa=\max\{\kappa_1,\kappa_2\}$ satisfies $\kappa<1$ in every case. We start with a test problem whose exact solution has a singularity at $t=0$. The availability of the exact solution will directly show that how accurately q-HWCM reflects the theoretical convergence near the point of singularity in practice.

\begin{example}\label{example_1}
Consider the following Caputo-type \(q\)-fractional coupled system defined on \(0<\zeta\leq 1\):

\begin{equation}\label{eq:ex1-system}
\begin{aligned}
{}^{C}D_q^{\alpha_1}\omega_1(\zeta)
&+\theta_1\zeta^{-\gamma_1}{}^{C}D_q^{\beta_1}\omega_1(\zeta)=f_1\bigl(\zeta,\omega_1(\zeta),\omega_2(\zeta)\bigr)+\lambda_1 E_{q,\alpha_1}\bigl(\omega_1(\zeta)\bigr),\\[3pt]
{}^{C}D_q^{\alpha_2}\omega_2(\zeta)
&+\theta_2\zeta^{-\gamma_2}{}^{C}D_q^{\beta_2}\omega_2(\zeta)=f_2\bigl(\zeta,\omega_1(\zeta),\omega_2(\zeta)\bigr)+\lambda_2 E_{q,\alpha_2}\bigl(\omega_2(\zeta)\bigr).
\end{aligned}
\end{equation}

subject to
\begin{equation*}
	\omega_1(0)=a_1=0.1,
	\qquad
	\omega_2(0)=a_2=0.05.
\end{equation*}
\end{example}
\noindent We did the experiments for different choices of the fractional order $\alpha_1,~\alpha_2$, whereas the other parameters are selected as: $q=0.7$, $\beta_1=0.3$, $\gamma_1=0.2$, $\mu_1=\alpha_1-\beta_1=0.5$, $\beta_2=0.2$, $\gamma_2=0.1$, $\mu_2=\alpha_2-\beta_2=0.5$, and $\theta_1=\theta_2=\lambda_1=\lambda_2=\varepsilon_1=\varepsilon_2=10^{-3}$. To construct the exact solution profile, we take the help of the following polynomial:
\[P(\zeta)=\sum_{m=1}^{12}p_m\zeta^m,\]
where the coefficient vector is
\[
\begin{aligned}
(p_1,\ldots,p_{12})=(&-0.241,-0.198,0.680,0.020,-0.181,-0.098,-0.883,0.014,1,-0.005,0.155,-0.397).
\end{aligned}
\]
The exact state variables are then defined through \(q\)-fractional integration:
\[\omega_1(\zeta)=a_1+J_q^{\alpha_1}g_1(\zeta),\qquad\omega_2(\zeta)=a_2+J_q^{\alpha_2}g_2(\zeta),\]
where $g_1(\zeta)=0.02+0.005P(\zeta),~g_2(\zeta)=0.015+0.005P(\zeta)$. Let $a_{1,0}=0.02, a_{2,0}=0.015$, and $a_{1,m}=a_{2,m}=0.005p_m,~ 1\leq m\leq12$. With the help of $J_q^\nu \zeta^m=\frac{\Gamma_q(m+1)}{\Gamma_q(m+\nu+1)}\zeta^{m+\nu}$, we can consider the exact solution of Example \ref{example_1} as follows
\begin{equation}\label{eq:ex1-explicit-solutions}
\begin{aligned}
\omega_i(\zeta)
={}&a_i+\sum_{m=0}^{12}a_{i,m}\times\frac{\Gamma_q(m+1)}
{\Gamma_q(m+\alpha_i+1)}\zeta^{m+\alpha_i},
\quad i=1,2.
\end{aligned}
\end{equation}
Note that it has a time singularity at $t-0$. For this choice, we have
\[
\begin{aligned}
v_i(\zeta)
&:={}^{C}D_q^{\beta_i}\omega_i(\zeta)=J_q^{\mu_i}g_i(\zeta)=\sum_{m=0}^{12}a_{i,m}
\frac{\Gamma_q(m+1)}
{\Gamma_q(m+\mu_i+1)}\quad\times\zeta^{m+\mu_i}.
\end{aligned}
\]
The nonlinear source functions are manufactured as

\begin{equation}\label{eq:ex1-source-f1}
\begin{aligned}
f_1(\zeta,x_1,x_2)
={}&g_1(\zeta)+\theta_1\zeta^{-\gamma_1}v_1(\zeta)-\lambda_1E_{q,\alpha_1}\bigl(\omega_1(\zeta)\bigr)+\varepsilon_1\bigl(x_2-\omega_2(\zeta)\bigr).
\end{aligned}
\end{equation}

\begin{equation}\label{eq:ex1-source-f2}
\begin{aligned}
f_2(\zeta,x_1,x_2)
={}&g_2(\zeta)+\theta_2\zeta^{-\gamma_2}v_2(\zeta)-\lambda_2E_{q,\alpha_2}\bigl(\omega_2(\zeta)\bigr)+\varepsilon_2\bigl(x_1-\omega_1(\zeta)\bigr).
\end{aligned}
\end{equation}
Note that the above example becomes coupled as $f_1$ and $f_2$ involve $x_1$ and $x_2$ respectively. Substituting \(x_i=\omega_i(\zeta)\) into \eqref{eq:ex1-source-f1}--\eqref{eq:ex1-source-f2} makes the coupling correction terms vanish. Consequently,

\[
\begin{aligned}
f_i\bigl(\zeta,\omega_1(\zeta),\omega_2(\zeta)\bigr)
&+\lambda_i E_{q,\alpha_i}\bigl(\omega_i(\zeta)\bigr)=g_i(\zeta)+\theta_i\zeta^{-\gamma_i}v_i(\zeta),\hspace{4em}i=1,2.
\end{aligned}
\]

Since ${}^{C}D_q^{\alpha_i}\omega_i=g_i$ and ${}^{C}D_q^{\beta_i}\omega_i=v_i,$ the functions in \eqref{eq:ex1-explicit-solutions} satisfy \eqref{eq:ex1-system} exactly. Moreover, \(J_q^{\alpha_i}g_i(0)=0\), and hence \(\omega_i(0)=a_i\) for $i=1,2$. 

\begin{figure}[ht]
\centerline{%
\begin{tabular}{cc}
\includegraphics[height=5cm]{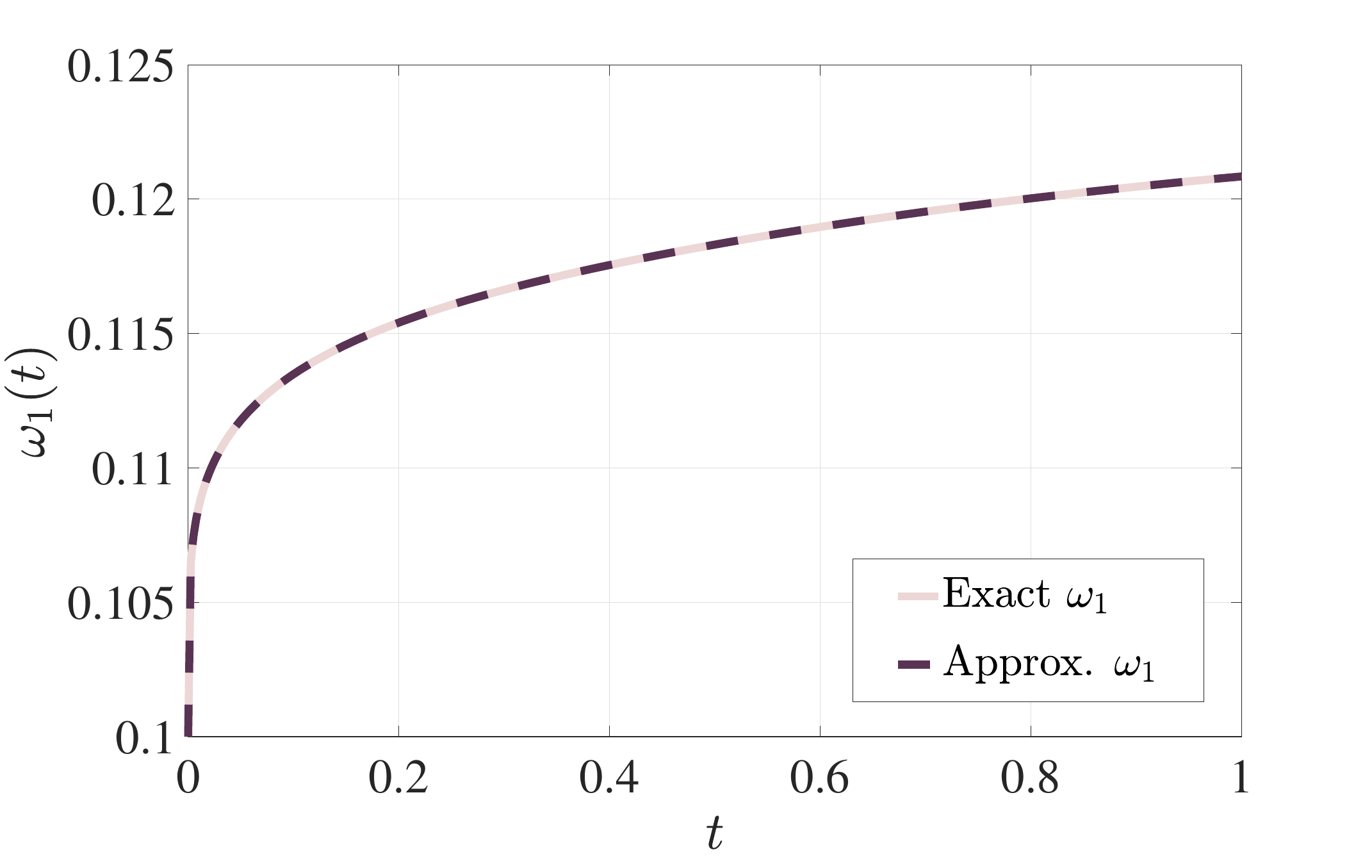} & \includegraphics[height=5cm]{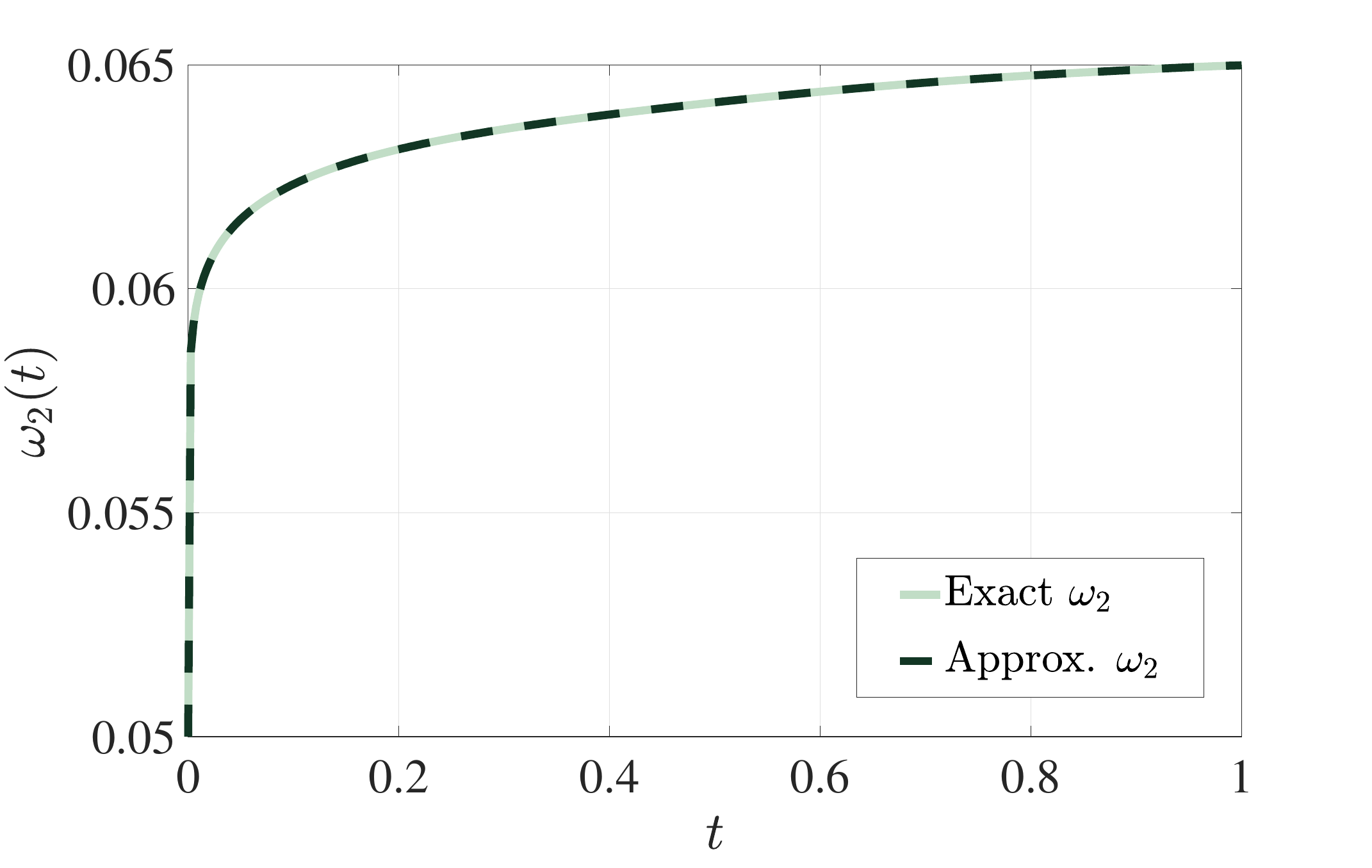}\\
(a) & (b)
\end{tabular}}
\caption{{\label{Ex1_Fig1} \it \small{Initial layer behavior in the numerical approximation for $\alpha_1=0.2,\alpha_2=0.1$ for Example \ref{example_1}.}}}
\end{figure}

\begin{figure}[ht]
\centerline{%
\begin{tabular}{cc}
\includegraphics[height=5cm]{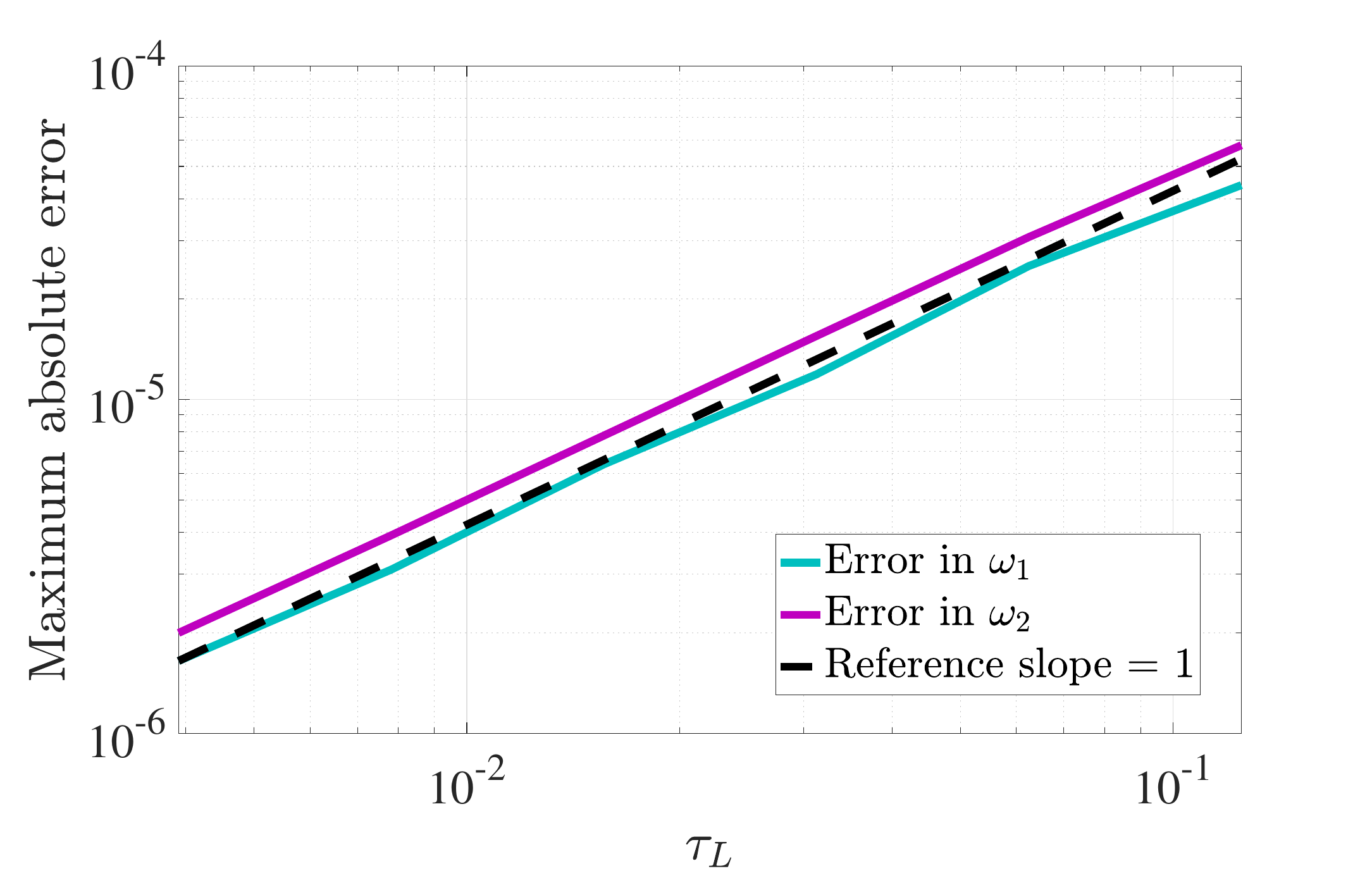} & \includegraphics[height=5cm]{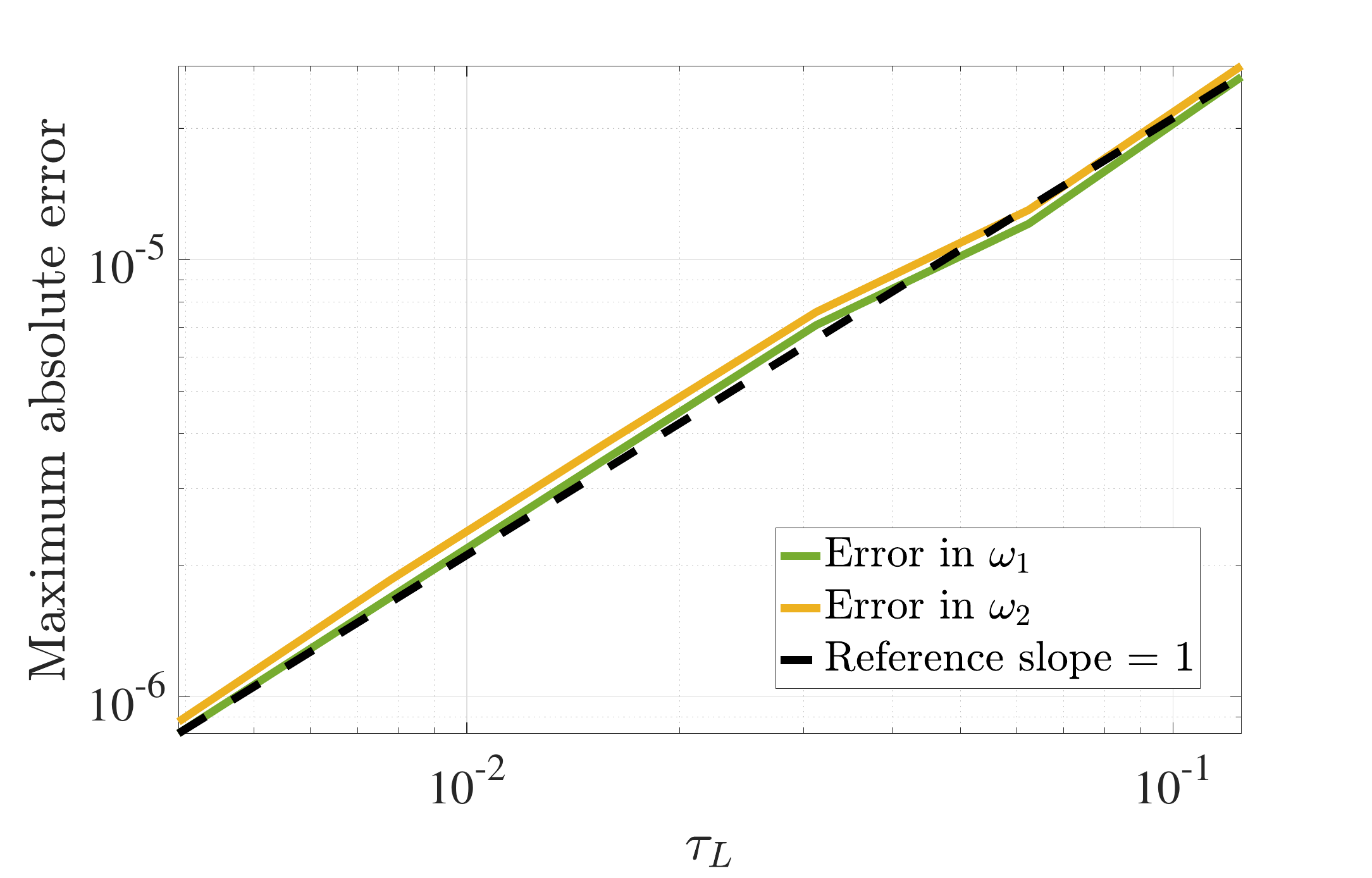}\\
(a) $\alpha_1=0.2,\alpha_2=0.1$ & (b) $\alpha_1=0.8,\alpha_2=0.7$
\end{tabular}}
\caption{{\label{Ex1_Fig2} \it \small{Convergence behavior of numerical approximations for Example \ref{example_1}.}}}
\end{figure}

\begin{figure}[ht]
\centerline{%
\begin{tabular}{cc}
\includegraphics[height=5cm]{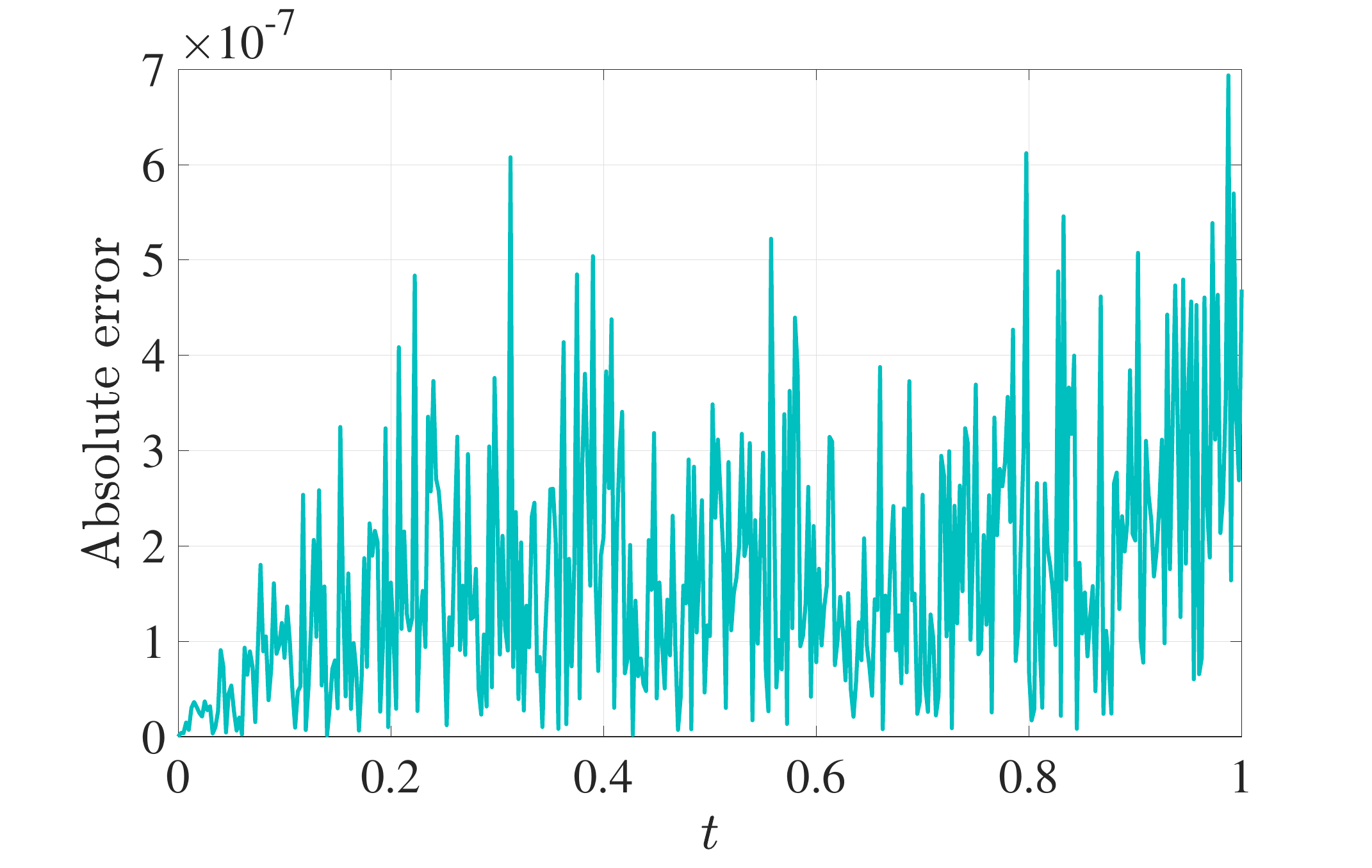} & \includegraphics[height=5cm]{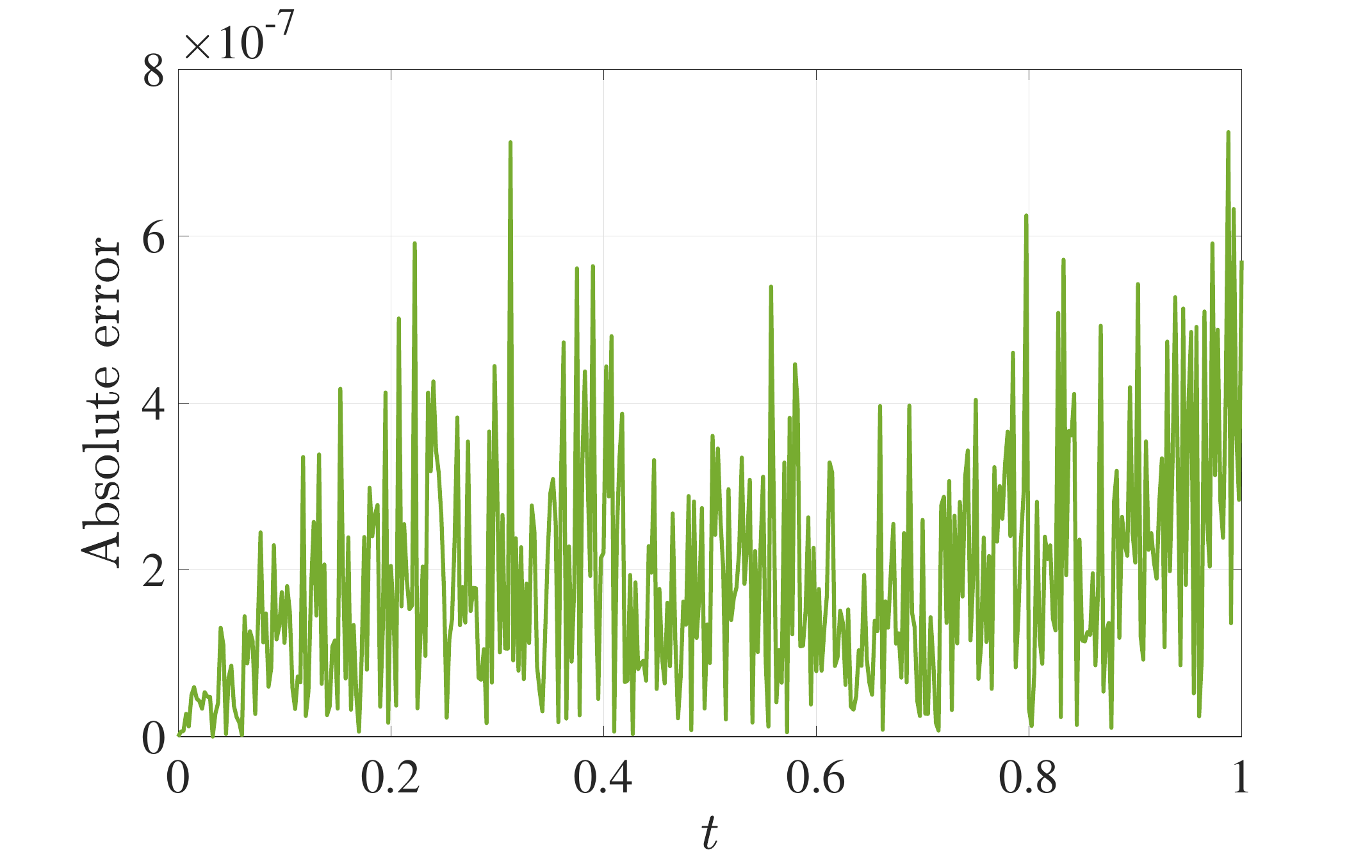}\\
(a) $\omega_1$ & (b) $\omega_2$
\end{tabular}}
\caption{{\label{Ex1_Fig3} \it \small{Absolute error of \(\omega_1\) and \(\omega_2\) at \(L=7\) $\alpha_1=0.8,\alpha_2=0.7$ for Example \ref{example_1}.}}}
\end{figure}

Figures \ref{Ex1_Fig1}–\ref{Ex1_Fig3} evaluate the numerical performance of the proposed Quantum Haar-Wavelet Collocation Method (q-HWCM) for Example \ref{example_1}, focusing on its initial-layer behavior, convergence properties, and pointwise absolute error distribution. Specifically, Figure \ref{Ex1_Fig1} contrasts the q-HWCM approximations with the exact solutions for the components $\omega_1(t)$ and $\omega_2(t)$ at $\alpha_1=0.2$ and $\alpha_2=0.1$. Across the entire domain $0\leq t\leq1$, the numerical and exact curves are visually indistinguishable, reflecting high fidelity in capturing the underlying fractional dynamics. Notably, both components undergo a sharp transition within a narrow vicinity of $t=0$ before evolving smoothly as $t$ increases—a classic initial-layer feature driven by fractional-order operators, especially when the fractional orders are small. Resolving this steep initial gradient is a well-known numerical challenge; the tight agreement here confirms that combining the localized, multiscale Haar-wavelet basis with a quantum-calculus formulation enables q-HWCM to capture rapid initial variations without spurious oscillations or degraded accuracy. Figure \ref{Ex1_Fig2} illustrates consistent convergence under grid refinement for both $(\alpha_1,\alpha_2)=(0.2,0.1)$ and $(0.8,0.7)$. In both cases, the error profiles align closely with the unit-slope reference line, confirming a roughly first-order convergence rate in $h_L$ across distinct fractional regimes. Finally, Figure \ref{Ex1_Fig3} displays the pointwise absolute errors for $L=7$ and $(\alpha_1,\alpha_2)=(0.8,0.7)$, showing that errors stay on the order of $10^{-7}$ across $0\leq t\leq1$ without noticeable amplification near $t=0$. Taken together, these plots verify that q-HWCM stably and accurately captures fractional initial layers while maintaining systematic convergence.Tables \ref{Ex1_Table1}–\ref{Ex1_Table4} further quantify the accuracy of the proposed q-HWCM. \\

Across all numerical tests, solutions and maximum errors are evaluated on a fixed uniform grid of $2049$ points over the computational domain, independent of the wavelet resolution level $L$. Table \ref{Ex1_Table1} demonstrates close agreement between the exact and numerical solutions for $\alpha_1=0.2$, $\alpha_2=0.1$, and $L=7$, with absolute errors for $\omega_1$ and $\omega_2$ staying below $10^{-6}$ at every tabulated point. Tables \ref{Ex1_Table2} and \ref{Ex1_Table3} detail the maximum errors $E_{i,L}$, their spatial locations, and the empirical convergence rates across two fractional-order pairs. For $(\alpha_1,\alpha_2)=(0.2,0.1)$, increasing $L$ from $2$ to $7$ reduces the maximum errors from $5.7854\times10^{-5}$ to $1.9997\times10^{-6}$ for $\omega_1$, and from $4.3944\times10^{-5}$ to $1.6479\times10^{-6}$ for $\omega_2$, with empirical rates remaining near unity. A comparable pattern emerges for $(\alpha_1,\alpha_2)=(0.8,0.7)$, where maximum errors drop to $8.2498\times10^{-7}$ and $8.7728\times10^{-7}$ at $L=7$ as the convergence rates approach one. These numerical findings align with the theoretical first-order error bound $\lVert\omega_i-\omega_{i,L}\rVert_\infty=\mathcal{O}(h_L)=\mathcal{O}(2^{-(L+1)})$ established in Theorem \ref{thm:convergence}, which holds even when $g_i\notin W^{1,\infty}(0,1)$. As shown in Table \ref{Ex1_Table4}, these results remain invariant even when the collocation grid is refined to cluster additional points near the origin's singularity. Ultimately, the tabular data reinforce the convergence trends in Figure \ref{Ex1_Fig2}, verifying that q-HWCM achieves reliable, first-order convergence across varying fractional orders while preventing error localization strictly within the initial layer.

\begin{table}[ht]
\centering
\caption{Exact and $q$-HWCM solutions for Example \ref{example_1}, where $\alpha_1=0.2,\alpha_2=0.1$, and $L=7$.}
\label{Ex1_Table1}
\footnotesize
\setlength{\tabcolsep}{3pt}
\begin{tabular}{c c c c c c c}
\toprule
$t$
& Exact $\omega_1$
& Approx. $\omega_1$
& Error
& Exact $\omega_2$
& Approx. $\omega_2$
& Error
\\
\midrule
0.1
& 0.11647515720987
& 0.11647562455513
& $4.6735\times10^{-7}$
& 0.06009558177800
& 0.06009598132193
& $3.9954\times10^{-7}$
\\
0.3
& 0.11818648996320
& 0.11818742633085
& $9.3637\times10^{-7}$
& 0.06240355342465
& 0.06240434417039
& $7.9075\times10^{-7}$
\\
0.5
& 0.11902656013408
& 0.11902620046245
& $3.5967\times10^{-7}$
& 0.06363084126363
& 0.06363061860852
& $2.2266\times10^{-7}$
\\
0.7
& 0.11963320707609
& 0.11963290494146
& $3.0213\times10^{-7}$
& 0.06453201218050
& 0.06453172246418
& $2.8972\times10^{-7}$
\\
0.9
& 0.12004862757311
& 0.12004842172098
& $2.0585\times10^{-7}$
& 0.06520106291758
& 0.06520085520328
& $2.0771\times10^{-7}$
\\
\bottomrule
\end{tabular}
\end{table}

\begin{table}[ht]
\centering
\caption{Maximum errors, their locations, and the order of convergences for Example \ref{example_1}, with $\alpha_1=0.2,\alpha_2=0.1$.}
\label{Ex1_Table2}
\small
\begin{tabular}{c c c c c c c}
\toprule
$L$
& $E_{1,L}$
& $t_{\max,1}$
& $p_{1,L}$
& $E_{2,L}$
& $t_{\max,2}$
& $p_{2,L}$
\\
\midrule
2
& $5.7853793\times10^{-5}$
& 0.1245117
& 0.9144
& $4.3943857\times10^{-5}$
& 0.2495117
& 0.8061
\\
3
& $3.0695327\times10^{-5}$
& 0.1875000
& 0.9856
& $2.5132761\times10^{-5}$
& 0.1875000
& 1.0788
\\
4
& $1.5501410\times10^{-5}$
& 0.0937500
& 0.9910
& $1.1898442\times10^{-5}$
& 0.1245117
& 0.8954
\\
5
& $7.7989834\times10^{-6}$
& 0.1406250
& 0.9927
& $6.3967598\times10^{-6}$
& 0.2031250
& 1.0479
\\
6
& $3.9193747\times10^{-6}$
& 0.1015625
& 0.9708
& $3.0939766\times10^{-6}$
& 0.2343750
& 0.9088
\\
7
& $1.9997132\times10^{-6}$
& 0.1679688
& ---
& $1.6479075\times10^{-6}$
& 0.1953125
& ---
\\
\bottomrule
\end{tabular}
\end{table}

\begin{table}[ht]
\centering
\caption{Maximum errors, their locations, and the order of convergences for Example \ref{example_1}, where $\alpha_1=0.8,\alpha_2=0.7$.}
\label{Ex1_Table3}
\small
\begin{tabular}{c c c c c c c}
\toprule
$L$
& $E_{1,L}$
& $t_{\max,1}$
& $p_{1,L}$
& $E_{2,L}$
& $t_{\max,2}$
& $p_{2,L}$
\\
\midrule
2
& $2.6231863\times10^{-5}$
& 1.0000000
& 1.1158
& $2.7806229\times10^{-5}$
& 1.0000000
& 1.0948
\\
3
& $1.2104461\times10^{-5}$
& 1.0000000
& 0.7708
& $1.3019106\times10^{-5}$
& 1.0000000
& 0.7766
\\
4
& $7.0943390\times10^{-6}$
& 1.0000000
& 1.0285
& $7.5998330\times10^{-6}$
& 1.0000000
& 1.0069
\\
5
& $3.4776898\times10^{-6}$
& 1.0000000
& 1.0381
& $3.7818101\times10^{-6}$
& 1.0000000
& 1.0269
\\
6
& $1.6934685\times10^{-6}$
& 1.0000000
& 1.0376
& $1.8560375\times10^{-6}$
& 1.0000000
& 1.0811
\\
7
& $8.2497605\times10^{-7}$
& 0.9765625
& ---
& $8.7728394\times10^{-7}$
& 0.9765625
& ---
\\
\bottomrule
\end{tabular}
\end{table}

\begin{table}[ht]
\centering
\caption{Maximum errors and observed convergence rates for Example \ref{example_1}, where $\alpha_1=0.1,\alpha_2=0.2$.}
\label{Ex1_Table4}
\small
\begin{tabular}{c c c c c}
\toprule
$L$
& $E_{1,L}$
& $p_{1,L}$
& $E_{2,L}$
& $p_{2,L}$
\\
\midrule
 8 & $1.01909541\times10^{-6}$     &0.9919 & $8.31293862\times10^{-7}$     &0.9999
\\
9  & $5.12402046\times10^{-7}$    & 0.9829 & $4.15668455\times10^{-7}$  &   0.9726
\\
10  & $2.59255235\times10^{-7}$   &  1.0196 & $2.11812663\times10^{-7}$   &  1.0214
\\
11  &  $1.27875921\times10^{-7}$   &  0.9846 & $1.04345936\times10^{-7} $  &  0.9661
\\
12  &  $6.46232251\times10^{-8}$   &  --- & $5.34137137\times10^{-8}$   &  ---
\\
\bottomrule
\end{tabular}
\end{table}

The numerical results for Example \ref{example_1} demonstrate that the proposed q-HWCM accurately captures the fractional initial-layer behavior and exhibits consistent, approximately first-order convergence for different fractional-order parameters. We next consider Example \ref{example_2} to further assess the accuracy and applicability of the proposed method for a different set of problem parameters.
\begin{example}\label{example_2}
We consider the coupled Caputo-type \(q\)-fractional system

\begin{equation*}
\begin{aligned}
{}^{C}D_q^{\alpha_1}\omega_1(\zeta)
&+\theta_1\zeta^{-\gamma_1}{}^{C}D_q^{\beta_1}\omega_1(\zeta)=f_1\bigl(\zeta,\omega_1(\zeta),\omega_2(\zeta)\bigr)\quad+\lambda_1E_{q,\alpha_1}\bigl(\omega_1(\zeta)\bigr),\\[3pt]
{}^{C}D_q^{\alpha_2}\omega_2(\zeta)
&+\theta_2\zeta^{-\gamma_2}{}^{C}D_q^{\beta_2}\omega_2(\zeta)=f_2\bigl(\zeta,\omega_1(\zeta),\omega_2(\zeta)\bigr)\quad+\lambda_2E_{q,\alpha_2}\bigl(\omega_2(\zeta)\bigr).
\end{aligned}
\end{equation*}

for \(0<\zeta\leq1\), subject to
\begin{equation*}
		\omega_1(0)=a_1=0.08,
		\qquad
		\omega_2(0)=a_2=0.04.
\end{equation*}
\end{example}	
\noindent The parameters are fixed as $q=0.7$, $\beta_1=0.25$, $\mu_1=\alpha_1-\beta_1$, $\gamma_1=0.15$, $\beta_2=0.35$, $\mu_2=\alpha_2-\beta_2$, $\gamma_2=0.20$, $\theta_1=1.2\times10^{-3}$, $\theta_2=1.5\times10^{-3}$, $\lambda_1=10^{-3}$, $\lambda_2=8\times10^{-4}$, $\varepsilon_1=8\times10^{-4}$, $\varepsilon_2=7\times10^{-4}$, while the fractional orders \(\alpha_1\) and \(\alpha_2\) are varied across the different numerical assessments. Based on the polynomial $P(\zeta)$ as defined for Example \ref{example_1}, one can obtain the exact solution for Example \ref{example_2} as 

\begin{equation*}
\begin{aligned}
\omega_i(\zeta)
={}&a_i+\sum_{m=0}^{12}a_{i,m}\times\frac{\Gamma_q(m+1)}
{\Gamma_q(m+\alpha_i+1)}\zeta^{m+\alpha_i},
\quad i=1,2.
\end{aligned}
\end{equation*}

where $a_{1,0}=0.018,a_{2,0}=0.012$, and $a_{1,m}=0.004p_m,a_{2,m}=-0.0035p_m,~ 1\leq m\leq12$. The lower-order derivatives can then be estimated as

\begin{equation*}
\begin{aligned}
v_i(\zeta)
&:={}^{C}D_q^{\beta_i}\omega_i(\zeta)=J_q^{\mu_i}g_i(\zeta)=\sum_{m=0}^{12}a_{i,m}
\frac{\Gamma_q(m+1)}
{\Gamma_q(m+\mu_i+1)}\times\zeta^{m+\mu_i},
\quad i=1,2.
\end{aligned}
\end{equation*}

The nonlinear source functions are defined as

\begin{equation*}
\begin{aligned}
f_1(\zeta,x_1,x_2)
={}&g_1(\zeta)+\theta_1\zeta^{-\gamma_1}v_1(\zeta)-\lambda_1E_{q,\alpha_1}\bigl(\omega_1(\zeta)\bigr)+\varepsilon_1\bigl[x_1^2-\omega_1(\zeta)^2+x_2-\omega_2(\zeta)\bigr],\\[3pt]
f_2(\zeta,x_1,x_2)
={}&g_2(\zeta)+\theta_2\zeta^{-\gamma_2}v_2(\zeta)-\lambda_2E_{q,\alpha_2}\bigl(\omega_2(\zeta)\bigr)+\varepsilon_2\bigl[\sin x_1-\sin\omega_1(\zeta)+x_2^2-\omega_2(\zeta)^2\bigr].
\end{aligned}
\end{equation*}

\begin{figure}[ht]
\centerline{%
\begin{tabular}{cc}
\includegraphics[height=4.7cm]{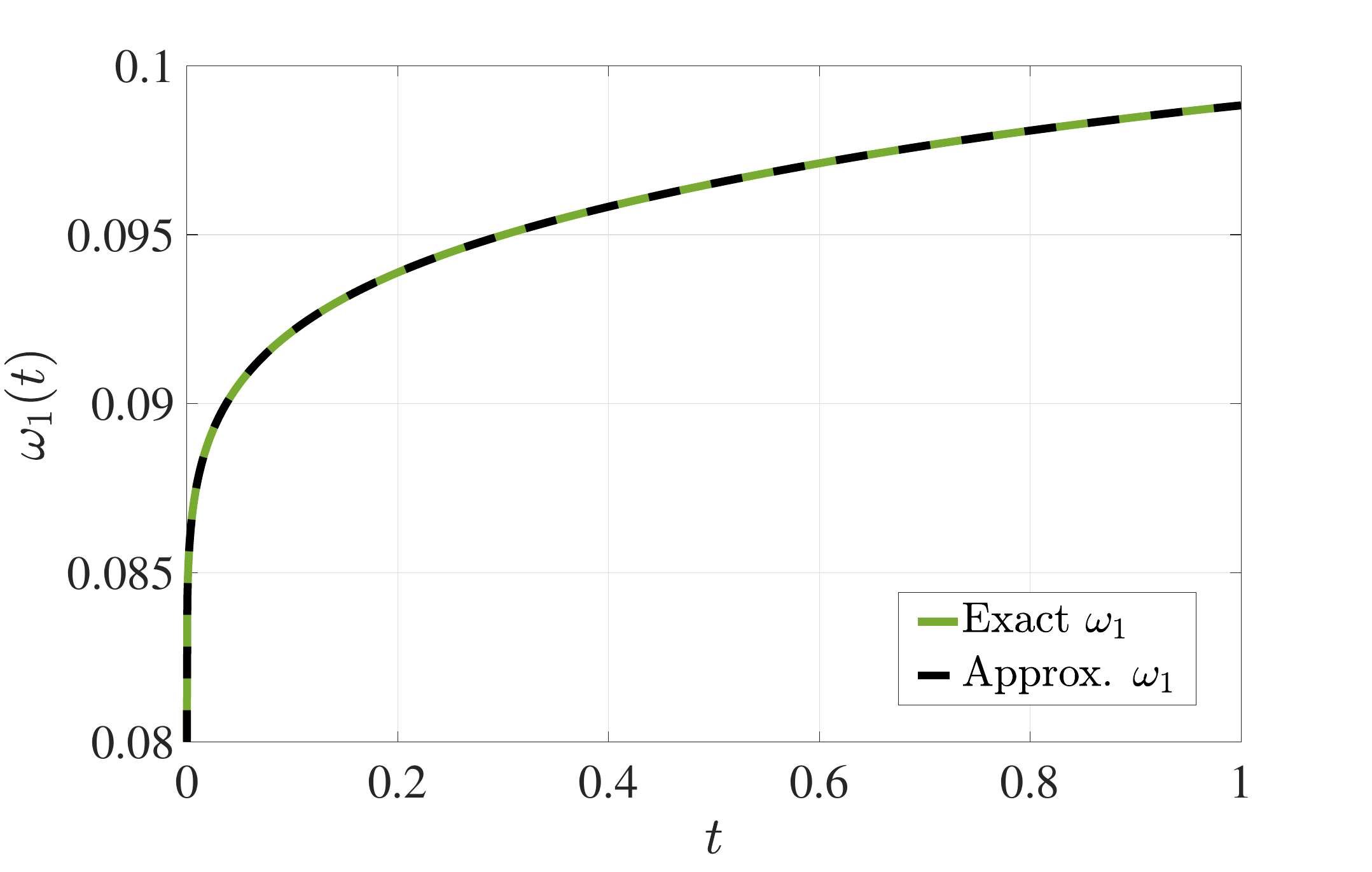} & \includegraphics[height=4.7cm]{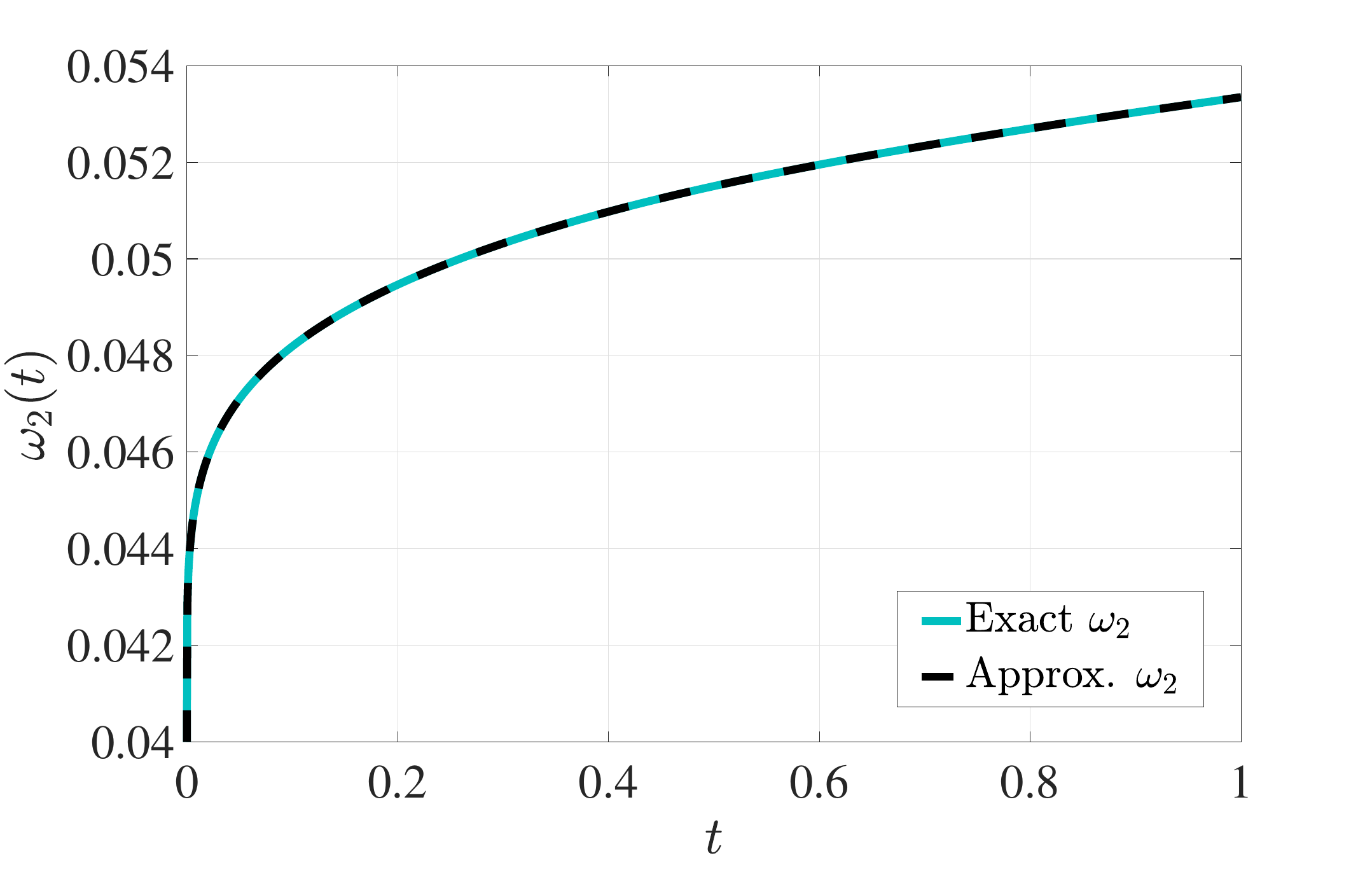}\\
(a) & (b)
\end{tabular}}
\caption{{\label{Ex2_Fig1} \it \small{Initial layer behavior with $\alpha_1=\alpha_2=0.2$ for Example \ref{example_2}.}}}
\end{figure}

\begin{figure}[ht]
\centerline{%
\begin{tabular}{cc}
\includegraphics[height=4.7cm]{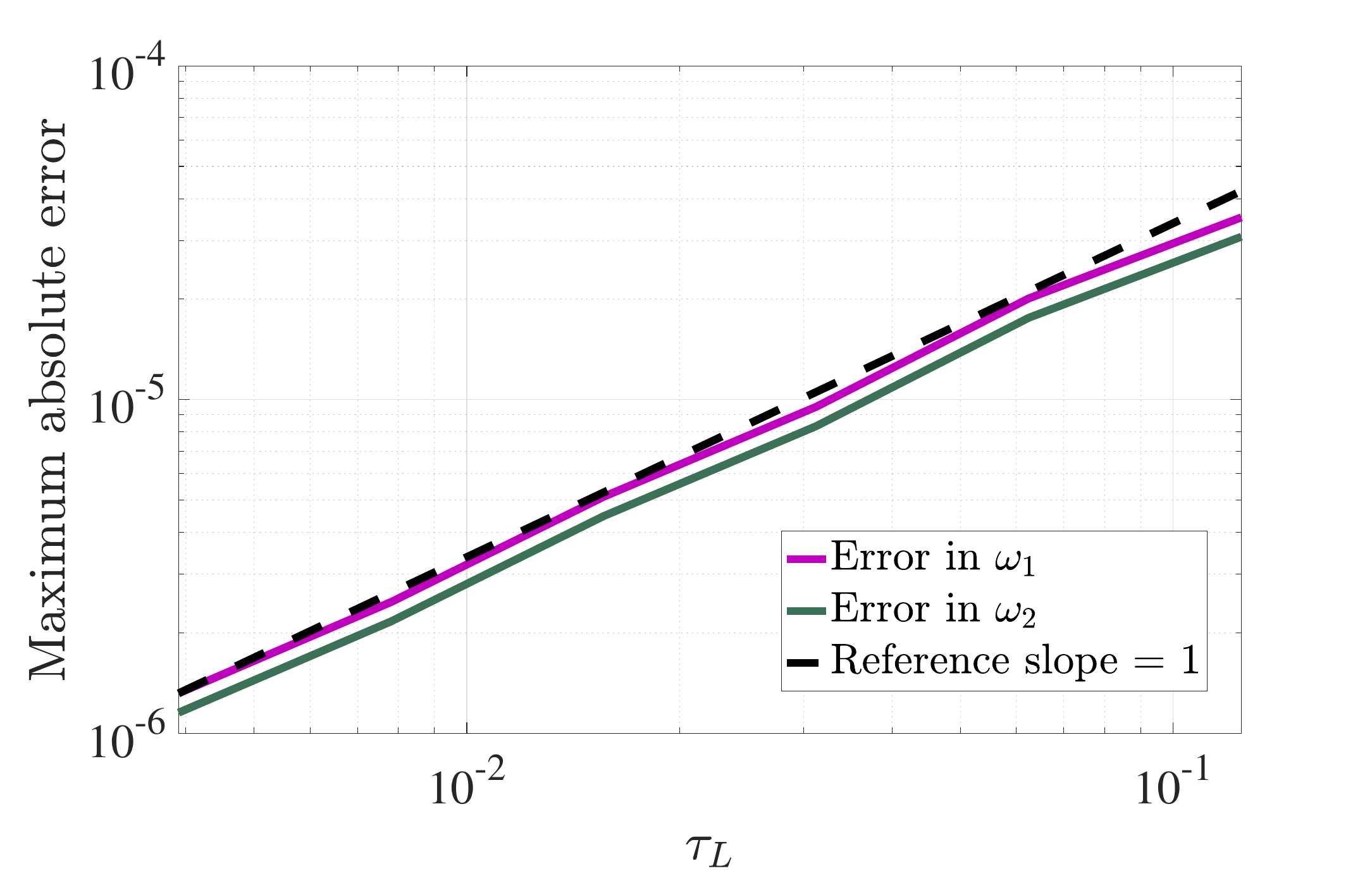} & \includegraphics[height=4.7cm]{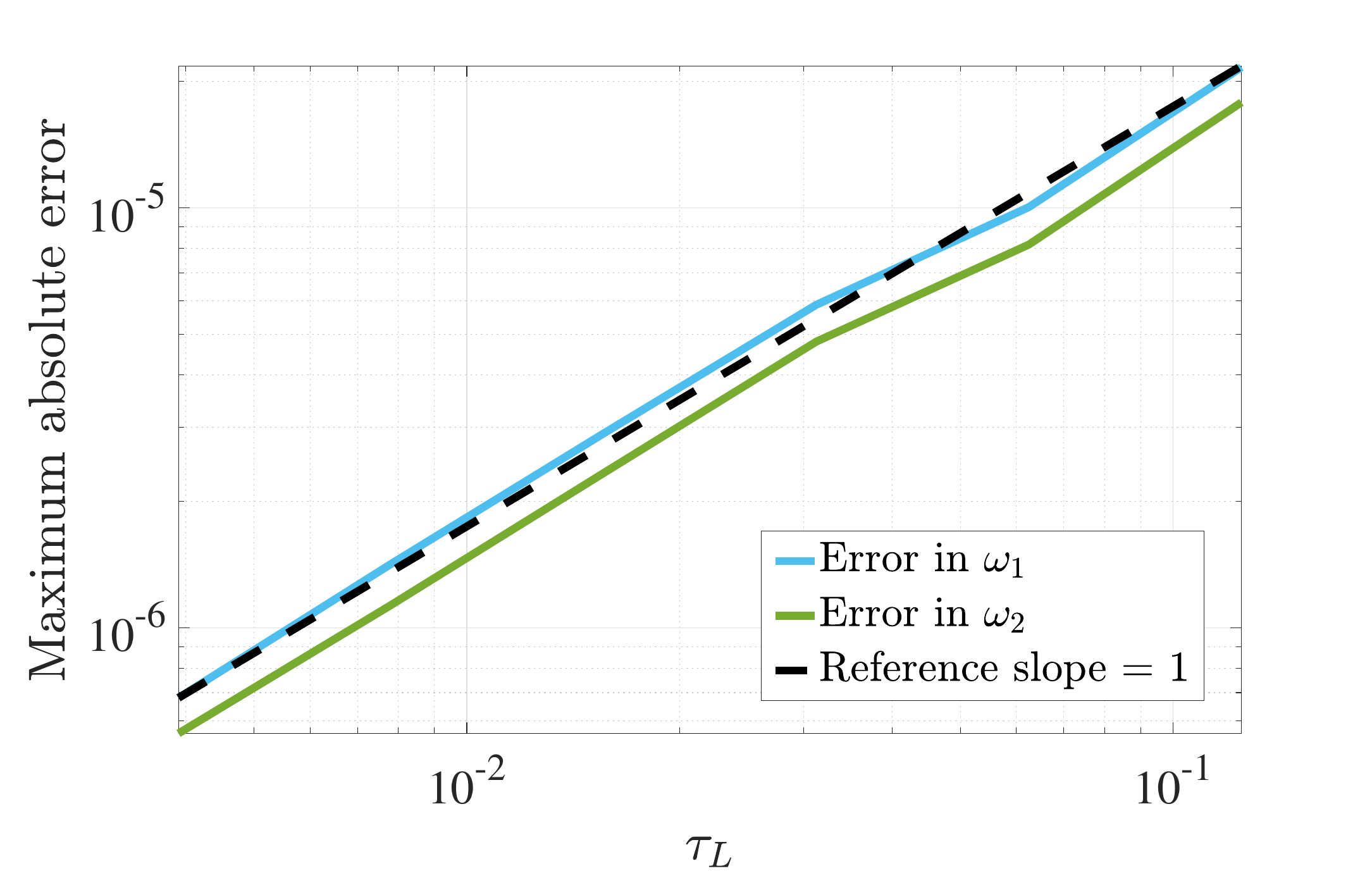}\\
(a) $\alpha_1=\alpha_2=0.2$ & (b) $\alpha_1=0.75,\alpha_2=0.85$
\end{tabular}}
\caption{{\label{Ex2_Fig2} \it \small{Convergence behavior as step length in space decreases, for Example \ref{example_2}.}}}
\end{figure}

\begin{figure}[ht]
\centerline{%
\begin{tabular}{cc}
\includegraphics[height=4.7cm]{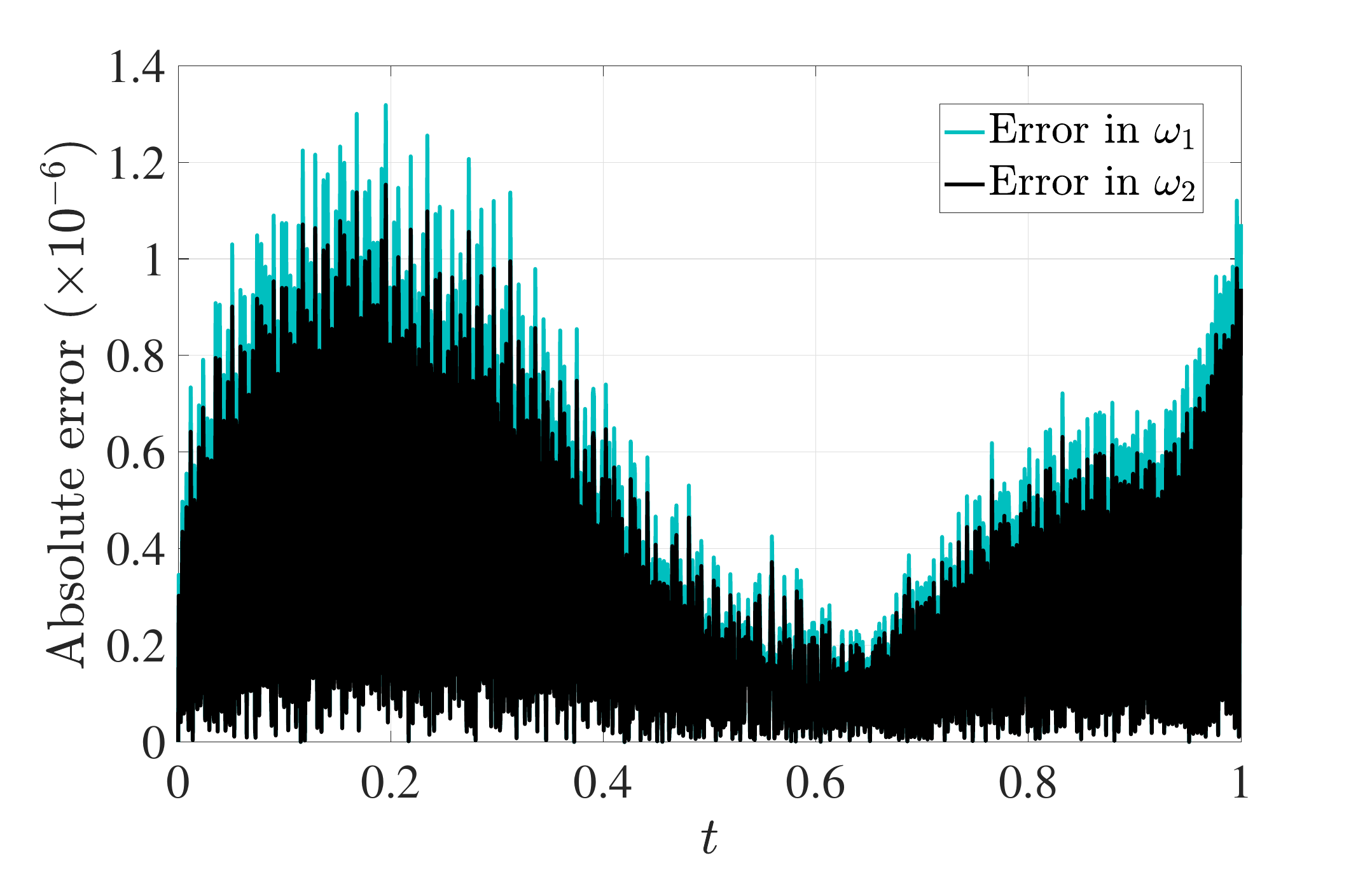} & \includegraphics[height=4.7cm]{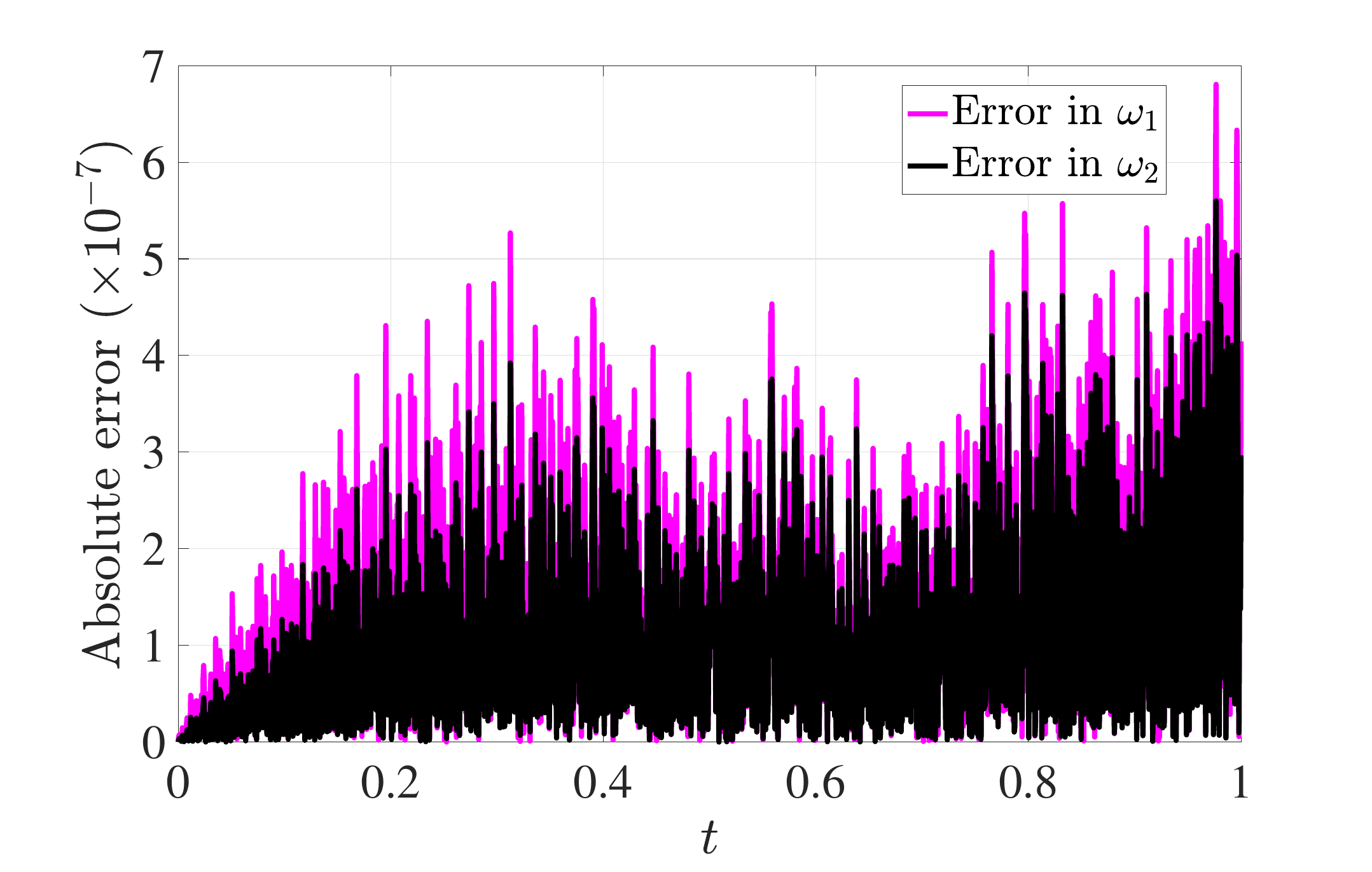}\\
(a) $\alpha_1=\alpha_2=0.2$ & (b) $\alpha_1=0.75,\alpha_2=0.85$
\end{tabular}}
\caption{{\label{Ex2_Fig3} \it \small{Absolute errors in \(\omega_1\) and \(\omega_2\) at \(L=7\) for Example \ref{example_2}.}}}
\end{figure}

Figures \ref{Ex2_Fig1}–\ref{Ex2_Fig3} together establish the accuracy of the q-HWCM scheme for Example \ref{example_2} across singular and mildly singular fractional-order regimes. In particular, Figure \ref{Ex2_Fig1} shows, for $\alpha_1=\alpha_2=0.2$, that both $\omega_1(t)$ and $\omega_2(t)$ exhibit a sharp initial-layer rise near $t=0$ followed by a slow concave increase to $t=1$ with the approximate (dashed) and exact (solid) curves indistinguishable throughout, confirming the method resolves the initial singularity rather than smoothing or oscillating around it. Figure \ref{Ex2_Fig2} quantifies convergence on a log–log plot of maximum absolute error versus $h_L$ for both $\alpha_1=\alpha_2=0.2$ and $\alpha_1=0.75,\alpha_2=0.85$, respectively. In both cases the $\omega_1$ and $\omega_2$ error curves decrease monotonically and run parallel to the reference slope-1 line across $h_L\in[10^{-3},10^{-1}]$, confirming first-order convergence that is not degraded by smaller $\alpha$, with slightly higher absolute error levels in the larger $\alpha$ due to larger solution magnitudes away from $t=0$. See the theoretical analysis given in Theorem \ref{thm:convergence} and further discussion in Remark \ref{remark_Convergence}. Figure \ref{Ex2_Fig3} fixes resolution at $L=7$ and shows the pointwise absolute error in time corresponding to \(N=256\) for $\alpha_1=\alpha_2=0.2$ the error is a dense oscillatory band peaking near $t\approx0.15$–$0.2$ (up to $\sim1.2$–$1.3\times10^{-6}$), dipping near $t\approx0.5$–$0.6$, and rising again symmetrically toward $t=1$, reflecting initial-layer error dominance. For $\alpha_1=0.75,\alpha_2=0.85$, the error is an order of magnitude smaller in scale ($\times10^{-7}$) but grows quasi-monotonically toward $t=1$ (up to $\sim5$–$6\times10^{-7}$), reflecting a weaker initial layer but greater endpoint/accumulated error as $\alpha$ increases. 

From the above tables, one can observe that the locations of maximum error slowly shifts towards outside the point of singularity once the order of the fractional derivative increases. However the order of accuracy of the present approach remains unchanged irrespective of the order of the fractional derivative.

Table \ref{Ex2_Table1} provides the pointwise maximum errors $E_{1,L}$ and $E_{2,L}$ between the exact and q-HWCM approximations of the components $\omega_2$ and $\omega_2$ respectively. It shows that the order of the error remains unaltered on a fixed resolution level at any regions including the region of singularity. However if one increases the resolution parameter of the Haar wavelet (see Tables \ref{Ex2_Table2}-\ref{Ex2_Table3}), error gradually decreases by retaining the same order of accuracy. One can also verify the theoretical rate of accuracy $O(2^{-L})$ from these two tables, as the resolution parameter increases $L$. In addition, note that the location of the maximum errors $t_{\max,1}$ and $t_{\max,2}$, are concentrated nearby to the initial time, when the order of the fractional derivative is relatively small, i.e., $\alpha_1=\alpha_2=0.2$, compared to the mildly singular case, when $\alpha_1=0.75$ and $\alpha_2=0.85$. Irrespective of the singularity behavior of the solution, the order of accuracy remains unchanged for all cases.


\begin{table}[ht]
\centering
\caption{Exact and q-HWCM solutions for Example \ref{example_2}, where $\alpha_1=\alpha_2=0.2$, and \(L=7\).}
\label{Ex2_Table1}
\resizebox{\textwidth}{!}{%
\begin{tabular}{ccccccc}
\toprule
\(t\)
& Exact \(\omega_1\)
& Approx. \(\omega_1\)
& Error
& Exact \(\omega_2\)
& Approx. \(\omega_2\)
& Error\\
\midrule
0.1 & 0.09214430363101 & 0.09214462321376
& \(3.1958\times10^{-7}\)
& 0.04818748603771 & 0.04818720641429
& \(2.7962\times10^{-7}\)\\
0.3 & 0.09499027465354 & 0.09499090730264
& \(6.3265\times10^{-7}\)
& 0.05032038227958 & 0.05031982870146
& \(5.5358\times10^{-7}\)\\
0.5 & 0.09651719004535 & 0.09651701172481
& \(1.7832\times10^{-7}\)
& 0.05150534999465 & 0.05150550607294
& \(1.5608\times10^{-7}\)\\
0.7 & 0.09762881619838 & 0.09762858443314
& \(2.3177\times10^{-7}\)
& 0.05233961567608 & 0.05233981847091
& \(2.0279\times10^{-7}\)\\
0.9 & 0.09847350705350 & 0.09847334093322
& \(1.6612\times10^{-7}\)
& 0.05303171865560 & 0.05303186399813
& \(1.4534\times10^{-7}\)\\
\bottomrule
\end{tabular}}
\end{table}

\begin{table}[ht]
\centering
\caption{Maximum errors, their locations, and the order of convergences for Example \ref{example_2}, with $\alpha_1=\alpha_2=0.2$.}
\label{Ex2_Table2}
\begin{tabular}{ccccccc}
\toprule
\(L\)
& \(E_{1,L}\)
& \(t_{\max,1}\)
& \(p_{1,L}\)
& \(E_{2,L}\)
& \(t_{\max,2}\)
& \(p_{2,L}\)\\
\midrule
2 & \(3.5162641\times10^{-5}\) & 0.2495117 & 0.8064
& \(3.0768808\times10^{-5}\) & 0.2495117 & 0.8065\\
3 & \(2.0105708\times10^{-5}\) & 0.1875000 & 1.0787
& \(1.7592416\times10^{-5}\) & 0.1875000 & 1.0787\\
4 & \(9.5189719\times10^{-6}\) & 0.1245117 & 0.8954
& \(8.3291302\times10^{-6}\) & 0.1245117 & 0.8954\\
5 & \(5.1172964\times10^{-6}\) & 0.2031250 & 1.0479
& \(4.4776157\times10^{-6}\) & 0.2031250 & 1.0479\\
6 & \(2.4750870\times10^{-6}\) & 0.2343750 & 0.9088
& \(2.1656817\times10^{-6}\) & 0.2343750 & 0.9088\\
7 & \(1.3182753\times10^{-6}\) & 0.1953125 & ---
& \(1.1534810\times10^{-6}\) & 0.1953125 & ---\\
\bottomrule
\end{tabular}
\end{table}

\clearpage
\begin{table}[ht]
\centering
\caption{Maximum errors, their locations, and the order of convergences for Example \ref{example_2}, where $\alpha_1=0.75,\alpha_2=0.85$.}
\label{Ex2_Table3}
\begin{tabular}{ccccccc}
\toprule
\(L\)
& \(E_{1,L}\)
& \(t_{\max,1}\)
& \(p_{1,L}\)
& \(E_{2,L}\)
& \(t_{\max,2}\)
& \(p_{2,L}\)\\
\midrule
2 & \(2.1605898\times10^{-5}\) & 1.0000000 & 1.1061
& \(1.7826074\times10^{-5}\) & 1.0000000 & 1.1251\\
3 & \(1.0037052\times10^{-5}\) & 1.0000000 & 0.7731
& \(8.1724605\times10^{-6}\) & 1.0000000 & 0.7680\\
4 & \(5.8731300\times10^{-6}\) & 1.0000000 & 1.0174
& \(4.7991046\times10^{-6}\) & 1.0000000 & 1.0393\\
5 & \(2.9014052\times10^{-6}\) & 1.0000000 & 1.0325
& \(2.3349913\times10^{-6}\) & 1.0000000 & 1.0442\\
6 & \(1.4183810\times10^{-6}\) & 1.0000000 & 1.0595
& \(1.1322598\times10^{-6}\) & 1.0000000 & 1.0159\\
7 & \(6.8054010\times10^{-7}\) & 0.9765625 & ---
& \(5.5992860\times10^{-7}\) & 0.9765625 & ---\\
\bottomrule
\end{tabular}
\end{table}

\section{Conclusion}\label{Section_Conclusion}
We developed and analyzed a $q$-Haar wavelet collocation method for coupled Caputo-type $q$-fractional Emden--Fowler systems with $q$-Mittag--Leffler nonlinearities. Expanding the highest-order fractional derivatives in Haar wavelets and reconstructing the lower-order terms via Jackson $q$-fractional integration avoids direct differentiation across the initial-time singularity and reduces the problem to a midpoint collocation system solved by damped Newton iteration. Using the Banach contraction principle, we proved well-posedness and perturbation stability for both the continuous and discrete formulations, and established $\mathcal{O}(\tau_L^\sigma)$ convergence in a weighted supremum norm that remains controlled up to $t=0$. Numerical experiments confirm first-order convergence on uniform dyadic partitions across both strongly and mildly singular regimes without requiring graded meshes or parameter tuning. Future work will extend this framework to higher-order $q$-wavelet bases, adaptive local refinement, and theoretical error bounds when $g_i$ itself has limited regularity.

\section*{Declaration}

\subsection*{Conflict of interest}
The authors declared that they have no conflict of interest.

\subsection*{Funding statement}
\noindent The author, Sudarshan Santra, would like to acknowledge the support of the Anusandhan National Research Foundation (ANRF), Govt. of India, for providing the research grant ANRF/ECRG/2025/002426/PMS.

\subsection*{Ethical approval}
In this manuscript, all the authors have agreed on authorship and approved the manuscript with consent for submission.

\subsection*{Credit authorship contribution statement}
Arushi: Writing - original draft, Software, Methodology, Investigation, Formal analysis, Review and editing the manuscript; Pranay Goswami: Visualization, Validation, Methodology, Investigation, Formal analysis, Conceptualization; Pratibhamoy Das: Investigation, Formal analysis, Writing - review \& editing, Visualization, Validation; Sudarshan Santra: Writing - original draft, Visualization, Validation, Investigation, Formal analysis, Data curation, Review and editing the manuscript.

\subsection*{Data availability statements}
No data was used for the research described in the article.

\subsection*{Declaration of generative AI use}
During the preparation of this work, the authors used ChatGPT to correct and improve the English language in several parts of the manuscript. After using this tool, the authors reviewed and edited the content as needed and take full responsibility for the content of the publication.

\FloatBarrier

\end{document}